\documentclass[11pt]{article}
\usepackage{amssymb}
\usepackage{url}
\usepackage{hyperref}
\usepackage{graphicx}
\usepackage{color}

\def\myn{{\mathfrak{n}}}
\def\Row{{\mathrm{Row}}}
\def\half{\hbox{\small${1\over 2}$}}
\def\four{\hbox{\small${1\over 4}$}}

\def\RiskOpt{{\hbox{\rm RiskOpt}}}
\def\diag{{\hbox{\rm diag}}}
\def\mR{{{\mathfrak{R}}}}

\def\bX{{\mathbf{X}}}
\def\bP{{\mathbf{P}}}

\def\bU{{\mathbf{U}}}
\def\bY{{\mathbf{Y}}}

\def\s{{\hbox{\rm\scriptsize s}}}

\def\bZ{{\mathbf{Z}}}

\def\Argmin{\mathop{\hbox{\rm Argmin}}}

\newtheorem{proposition}{Proposition}[section]

\definecolor{MyDarkBlue}{rgb}{0,0.08,0.45}
\definecolor{MyViolet}{rgb}{0.45,0.08,0.95}
\definecolor{MyBrown}{rgb}{0.45,0.08,0}

\def\norm2to2{{\|\cdot\|_{2,2}}}
\def\Prob{\hbox{\rm Prob}}

\def\cl{{\hbox{\rm  cl}\,}}

\def\bE{{\mathbf{E}}}

\def\inter{\hbox{\rm  int}}

\def\Diag{\hbox{\rm  Diag}}
\def\Prob{\hbox{\rm  Prob}}

\def\Opt{\hbox{\rm Opt}}

\def\Conv{\hbox{\rm  Conv}}

\def\Tr{{\mathop{\hbox{\rm  Tr}}}}
\def\cA{{\cal A}}
\def\cB{{\cal B}}

\def\cE{{\cal E}}

\def\cH{{\cal H}}

\def\cM{{\cal M}}
\def\cN{{\cal N}}

\def\cQ{{\cal Q}}
\def\cR{{\cal R}}
\def\cS{{\cal S}}
\def\cT{{\cal T}}

\def\cX{{\cal X}}
\def\cY{{\cal Y}}
\def\cZ{{\cal Z}}
\def\SG{{\cal SG}}

\def\mDelta{{\mathbf \Delta}}

\def\Argmin{\mathop{\hbox{\rm  Argmin}}}

\def\abs{\mbox\hbox{\rm  abs}}

\def\bS{{\mathbf{S}}}
\def\R{{\mathbb R}}

\def\bH{{\mathbf{H}}}
\def\e{{\hbox{\rm e}}}
\def\abs{{\hbox{\rm abs}}}

\def\qed{\ \hfill$\square$\par\smallskip}

\def\mypict3{\epsfxsize=220pt\epsfysize=80pt\epsfbox}

\def\bR{{\mathbf{R}}}

\newtheorem{lemma}{Lemma}[section]
\def\cH{{\cal H}}
\def\Col{{\hbox{\rm Col}}}

\def\Risk{{\hbox{\rm Risk}}}

\newcommand{\aic}[2]{{\color{MyDarkBlue}~#2}}

\newcommand{\be}{\begin{eqnarray}}
\newcommand{\ee}[1]{\label{#1}\end{eqnarray}}
\newcommand{\nn}{\nonumber \\}
\newcommand{\ese}{\end{eqnarray*}}
\newcommand{\bse}{\begin{eqnarray*}}
\newcommand{\rf}[1]{~(\ref{#1})}
\newcommand{\wh}[1]{{\widehat{#1}}}

\title{On Polyhedral Estimation of Signals via Indirect Observations}
\author{
Anatoli Juditsky
\thanks{LJK, Universit\'e Grenoble Alpes, 700 Avenue Centrale 38401 Domaine Universitaire
de Saint-Martin-d'H\`{e}res, France,
{\tt anatoli.juditsky@univ-grenoble-alpes.fr}}
\and Arkadi Nemirovski
\thanks{Georgia Institute
 of Technology, Atlanta, Georgia
30332, USA, {\tt nemirovs@isye.gatech.edu}\newline
The first author was supported by the PGMO grant 2016-2032H. Research of
the second author was supported by NSF grant   CCF-1523768.}}
\date{}
\begin{document}
\maketitle
\begin{abstract}
We consider the problem of recovering linear image of unknown signal belonging to a given convex compact signal set from noisy observation of another linear image of the signal. We develop a simple generic efficiently computable {\sl non}linear in observations ``polyhedral'' estimate along with computation-friendly techniques for its design and risk analysis. We demonstrate that under favorable circumstances the resulting estimate is provably near-optimal in the minimax sense, the ``favorable circumstances'' being less restrictive than the weakest known so far assumptions ensuring near-optimality of estimates which are linear in observations.
\end{abstract}
\section{Introduction}
\subsection{Motivation}\label{plhdmotiv}
In this paper we consider the  estimation problem as follows:
 \begin{quote} {\em Given a noisy observation
\be
\omega=Ax+\xi\in\bR^m,
\ee{eq:mod1}
of linear image $Ax$ of an unknown signal $x$, we want to recover the image $w=Bx\in\bR^\nu$ of this signal. It is assumed that $x$ is known to belong to a given signal set -- a nonempty convex compact set $\cX\subset\bR^n$, and $A$ and $B$ are given $m\times n$ and $\nu\times n$ matrices; $\xi$ is the observation noise with distribution $P_x$ which may depend on $x$.  In addition, we are given a norm $\|\cdot\|$ on $\bR^\nu$ in which the estimation error is measured.}
\end{quote}
Estimation problem \rf{eq:mod1} is a (nonparametric) linear inverse problem. Popular approaches to solving statistical inverse problems include, among others, Singular Value Decomposition (SVD), see, e.g., \cite{johnstone1990speed,johnstone1991discretization,mair1996statistical},
Galerkin projection\cite{natterer1986mathematics}, Wavelet-Vaguelette Decompositions \cite{donoho1995nonlinear,abramovich1998wavelet,figueiredo2003algorithm,johnstone2004wavelet}, adaptive Galerkin algorithms \cite{cohen2004adaptive,hoffmann2008nonlinear}, as well as various iterative regularization techniques \cite{vogel2002computational,kaipio2006statistical}. When statistically analysed, those approaches usually assume a special structure of the problem, when matrix $A$ and set $\cX$ ``fit each other,'' e.g., there exists a sparse approximation of the set $\cX$ in a given basis/pair of bases, in which matrix $A$ is ``almost diagonal'' (see, e.g. \cite{donoho1995nonlinear,cohen2004adaptive} for detail). Under these assumptions, traditional results focus on estimation algorithms which are both numerically straightforward and statistically (asymptotically) optimal with closed form analytical description of estimates and corresponding risks.
In the situation considered in this paper, where we do not assume any specific structure of $\cX$ apart from convexity, compactness and ``computationally tractability,''\footnote{For a brief outline of computational tractability, see Section \ref{comptract} of the appendix.} and $A$, $B$ are ``general'' matrices of appropriate dimensions,
 we cannot expect deriving closed form expressions for estimates and risks. Instead, we adopt an alternative approach initiated in \cite{Don95} and further developed in \cite{JN2009,GJN2015,l2estimation,JudNem2018}. Within this {\em operational} approach both the estimate and its risk are yielded by efficient computation, usually via convex optimization, rather than by an explicit closed form analytical description; what we know in advance, in good cases,  is that the resulting risk, whether large or low, is nearly the best one achievable under the circumstances.

Note that a ``standard choice'' of statistical techniques which can be applied to solve \rf{eq:mod1} is between Maximum Likelihood and Linear Estimation. Maximal Likelihood Estimation (MLE)
\cite{fisher1912absolute,fisher1922mathematical,pratt1976fy} is a ``universal'' statistical tool; when the distribution $P$ of the noise is known, the application  to problem \rf{eq:mod1} is straightforward.
It is well known that MLE has excellent asymptotical properties in parametric models; however, it is also well known \cite{bahadur1958examples,le1990maximum,birge2006model}
that Maximum Likelihood approach has significant drawbacks. For instance, let us consider the simple situation of direct observation where $A=B=I_n$,
observations noise is Gaussian -- $\xi\sim \cN(0,I_n)$, $\|\cdot\|=\|\cdot\|_2$,  and
\[
\cX=\big\{x=[x_1;...;x_{n}]:\,|x_1|\leq n^{1/4},\,\|[x_2;...;x_{n}]\|_2+2n^{-1/4}|x_1|\leq 2\big\}.
\]
In this case,  for $n$ large enough, the squared risk of the MLE $\wh x^{ML}$ satisfies \cite{birge2006model}
\[
\sup_{x\in \cX}\bE
\{\|x-\wh x^{ML}\|_2^2\}\geq \mbox{$3\over 4$}\sqrt{n},
\]
while for the minimax risk it holds
\be
\inf_{\wh x}\sup_{x\in \cX}\bE
\{\|x-\wh x\|_2^2\}\leq 5,
\ee{eq:lingood}
with the upper bound \rf{eq:lingood} attained at the simple linear estimator $\wh x=[\omega_1;0;...;0]$.

Aside of Maximum Likelihood, very popular, due to their relative simplicity, estimates are the {\sl linear ones} -- those of the form $\omega\mapsto\widehat{w}(\omega)=G^T\omega$. Linear estimates have received much attention in the statistical literature  (cf. \cite{kuks1,kuks2,Rao1973,Pinsker1980,donoho1990minimax,efromovich1996sharp,Cris,wasserman2006all,Tsybakov} among many others)
When designing a linear estimate, the emphasis is on how to specify the matrix $G$ in order to  obtain the lowest possible maximal over $\cX$ estimation risk, which is then compared to the minimax risk (the infimum, taken w.r.t. all Borel estimates $\widehat{x}(\cdot)$, of the worst case, over signals from $\cX$, expected norm of the recovery error). ``Near optimality'' results for the case of indirect observations (where $A$ and $B$ are arbitrary) are  the subject of recent papers \cite{l2estimation,JudNem2018}, where it was shown that in the {\sl spectratopic case}, where $\cX$ and the unit ball $\cB_*$ of the norm conjugate to $\|\cdot\|$ are {\sl spectratopes}\footnote{see \cite{JudNem2018} or Section \ref{plhdcompat}.A below; as of now, a instructive example of spectratope is the  intersection of finite family of ellipsoids/elliptic cylinders with common center, or, to give  a more exotic example, the unit ball of the spectral norm in the space of matrices.},  a properly designed, via solving an explicit convex optimization problem, linear estimate is nearly optimal.
\par
What follows is motivated by the desire to build an alternative estimation scheme which works beyond the ellitopic/spectratopic case, where linear estimates can become ``heavily nonoptimal.''
\bigskip\par\noindent\textbf{Motivating example.} Consider the simple-looking problem of recovering $Bx=x$ in $\|\cdot\|_2$-norm from {\sl direct} observations $(Ax=x)$ corrupted by the standard Gaussian noise $\xi\sim\cN(0,\sigma^2 I)$, and let $\cX$ be the unit $\|\cdot\|_1$=ball:
$$
\cX=\Big\{x\in\bR^n: \sum_i|x_i|\leq1\Big\}.
$$
In this situation, building the optimal, in terms  of the worst-case over $x\in \cX$ expected $\|\cdot\|_2^2$-risk, linear estimate $\widehat{x}_H(\omega)=H^T\omega$ is extremely simple.
Indeed, one easily verifies that
the optimal $H$ is a scalar matrix $h I$, with the optimal $h$ being the minimizer of the univariate quadratic function $(1-h)^2+\sigma^2nh^2$. Therefore, in this case, the best achievable  with linear estimates expected $\|\cdot\|_2^2$-risk is
$$
\max_{x\in \cX}\bE\left\{\|\widehat{x}_H(\omega)-x\|_2^2\right\}=\min_h\left[(1-h)^2+\sigma^2nh^2\right]  ={n\sigma^2\over 1+n\sigma^2}.
$$
On the other hand, consider {\sl nonlinear} estimate as follows. Given  observation $\omega$, we specify the estimate $\widehat{x}(\omega)$ as an optimal solution to the optimization problem
\be
\Opt(\omega)=\min_{y\in\cX} \|y-\omega\|_\infty.
\ee{0Eq1}
Note that for every $\rho>0$ the probability for the true signal to satisfy $\|x-\omega\|_\infty\leq\rho\sigma$ (``event $\cE$'') is at least $1-\epsilon$ for $\epsilon=2n\exp\{-\rho^2/2\}$, and if this event occurs then
 both $x$ and $\widehat{x}$ belong to the box $\{y:\|y-\omega\|_\infty\leq\rho\sigma\}$, implying that $\|x-\widehat{x}\|_\infty \leq 2\rho\sigma$. Combining the latter bound with the constraint
 $\|x-\widehat{x}\|_2\leq\|x-\widehat{x}\|_1\leq 2$, since $x\in\cX$ and $\widehat{x}\in\cX$, we  obtain
 $$
 \|x-\widehat{x}\|_2\leq \sqrt{\|x-\widehat{x}\|_\infty\|x-\widehat{x}\|_1}\leq\left\{\begin{array}{ll}2\sqrt{\rho\sigma},&\omega\in\cE\\
 2,&\omega\not\in \cE\\
 \end{array}\right.,
 $$
whence
$$
\bE\left\{\|\widehat{x}-x\|_2^2\right\}\leq 4\rho\sigma+4\epsilon\leq 4\rho\sigma+8n\exp\{-\rho^2/2\}.\eqno{(*)}
$$
Assuming $\sigma\leq 2n/\sqrt{\e}$  and specifying $\rho$ as $\sqrt{2\ln(2n/\sigma)}$, we get $\rho\geq1$ and $2n\exp\{-\rho^2/2\}\leq\sigma$, implying that the right
hand side in $(*)$ is at most $8\rho\sigma$. Therefore, the nonlinear estimate $\widehat{x}(\omega)$ satisfies
$$
\max_{x\in \cX}\bE\left\{\|\widehat{x}(\omega)-x\|_2^2\right\}\leq 8\sqrt{\ln(2n/\sigma)}\sigma.
$$ When $n\sigma^2$ is of order of 1, the latter bound  is of order of
$\sigma\sqrt{\ln(1/\sigma)}$, while the best expected $\|\cdot\|_2^2$-risk attainable with linear estimates under the circumstances is of order of 1.
We conclude that when $\sigma$ is small and $n$ is large (specifically, is of order of $1/\sigma^2$),  the best linear estimate is {\sl far inferior} to the nonlinear estimate -- the ratio of the corresponding squared risks is as large as ${O(1)\over\sigma\sqrt{\ln(1/\sigma)}}$.

\subsection{Polyhedral estimate}
The construction of the nonlinear estimate $\widehat{x}$ we have built in the above example\footnote{In fact, this estimate is nearly optimal under the circumstances in a meaningful range of values of $n$ and $\sigma$.} admits a natural extension yielding what we call  {\em polyhedral estimate.}
The idea underlying polyhedral estimate is quite straightforward. Assuming for present that the observation noise is $\cN(0,\sigma^2I_n)$, observe that there is a spectrum of ``easy to estimate'' linear forms of signal $x$ underlying observation, namely the forms $g_h^Tx=h^TAx$ with $h\in\cH=\{h\in\bR^m:\|h\|_2=1\}$. Indeed, for a form of this type, the ``plug-in'' estimate $\widehat{g}_h(\omega)=h^T\omega$ is  an unbiased estimate of $g_h^Tx$ with $\cN(0,\sigma^2)$ recovery error. It follows that selecting somehow a {\sl contrast matrix} $H$ -- an $m\times M$ matrix with columns from $\cH$, the plug-in estimate $H^T\omega$ recovers well the vector $H^TAx$ in the uniform norm:
\begin{equation}\label{1Eq1}
\Prob_{\omega\sim\cN(Ax,\sigma^2I_m)}\left\{\|H^T\omega-H^TAx\|_\infty >\sigma\rho\right\}\leq 2M\exp\{-\rho^2/2\},\,\rho\geq0.
\end{equation}
As a result, given a ``reliability tolerance'' $\epsilon\ll1$ and setting $\rho=\sqrt{2\ln(2M/\epsilon)}$, the estimate $H^T\omega$ recovers the vector $H^TAx$, whatever be $x\in \bR^n$, within $\|\cdot\|_\infty$-accuracy $\sigma\rho$ and reliability  $1-\epsilon$. When our objective is to recover $w=Bx$, a natural way to combine this estimate with a priori information that $x\in \cX$ is to set
\begin{equation}\label{2Eq1}
\wh{x}^H(\omega)\in\Argmin_y\left\{\|H^T[Ay-\omega]\|_\infty:\,y\in\cX\right\},\;\;\widehat{w}^H(\omega)=B\wh{x}(\omega).
\end{equation}
Note that the estimate $\widehat{w}^H(\cdot)$ of $w$ we end up with is defined solely in terms of $H$ and the data $A,B,\cX$ of our
estimation problem, and that simple estimate (\ref{0Eq1}) is nothing but the polyhedral estimate stemming from the unit contrast matrix.
The rationale behind polyhedral estimation scheme is the desire to reduce
complex estimating problems to those of estimating linear forms.
To the best of our knowledge, the idea of polyhedral estimate  goes back to \cite{oldpaper}, see also \cite[Chapter 2]{saintflour}, where it was shown that when recovering smooth multivariate regression functions known to belong to Sobolev balls from their noisy observations taken along a regular grid $\Gamma$, a polyhedral estimate with ad hoc selected contrast matrix is near-optimal
in a wide range of smoothness characterizations and norms $\|\cdot\|$. Recently, the ideas underlying the results of \cite{oldpaper} have been taken up in the MIND estimator of \cite{grasmair2018variational}, then applied in the indirect observation setting in \cite{proksch2018multiscale} in the context of multiple testing.
\par
The goal of this paper is to investigate characteristics of the polyhedral estimate, with a particular emphasis on efficiently computable upper bounds for the risk of the estimate $\wh w^H(\cdot)$ and design of the {\em  contrast matrix $H$}  resulting in the (nearly) best upper risk bounds.
\par
The main body of the paper is organized as follows. We begin in Section \ref{sect1} with detailed formulation of the estimation problem (see Section \ref{theproblem}), and  present generic polyhedral estimate along with its risk analysis (Section \ref{polyhedest}). This analysis requires from the main ingredient of a polyhedral estimate -- the underlying {\sl contrast matrix} -- to be properly adjusted to the structure and the magnitude of observation noise. To allow for this adjustment, we restrict ourselves to three types of {\em observation schemes} {specifying the noise} structure, referred to as {\sl sub-Gaussian}, {\sl Discrete}, and {\sl Poisson} cases; description of these cases and of the restrictions they impose on the contrast matrices form the subject of Section \ref{setscH}. The subject of the subsequent sections is
tuning the polyhedral estimate to the structure of the estimation problem -- the design of the contrast matrix aimed at minimizing the risk of the associated polyhedral estimate. We propose two related approaches. The first of them, developed in Section \ref{secteffI}, mimics, in  a sense, the mechanism working in the above Motivating example. We show
 that this approach produces near-optimal estimates in some situations where linear estimates can be heavily nonoptimal. The second, completely different, approach to the design of contrast matrix is developed in Section \ref{secteffII}. Based on the notion of a ``cone compatible with convex set,'' this approach is inspired by the design of nearly optimal linear estimates in the case where the signal set and the unit ball of the norm conjugate for $\|\cdot\|$ are ellitopes/spectratopes. While this approach is {\sl not} restricted to the latter case, we show that in the spectratopic case, when the observation noise is zero mean sub-Gaussian, it results in polyhedral estimate which is provably nearly optimal in the minimax sense.
   \par
Technical proofs (which are longer than few lines) are relegated to the appendix.
\section{Problem of interest and generic polyhedral estimate}\label{sect1}
\subsection{The problem}\label{theproblem}
Suppose that we are given
\begin{itemize}
\item a nonempty computationally tractable convex compact {\sl signal set} $\cX\subset\bR^n$,
\item {\sl sensing matrix} $A\in\bR^{m\times n}$, {\sl decoding matrix} $B\in\bR^{\nu\times n}$, and a norm $\|\cdot\|$ on the space $\bR^\nu$,
\item a {\sl reliability tolerance} $\epsilon\in(0,1)$,
\item a random {\sl observation}
\begin{equation}\label{Eq-1}
\omega=Ax+\xi_x
\end{equation}
stemming from unknown {\sl signal} $x$ known to belong to $\cX$; here $\xi_x$ is a random variable with Borel probability distribution $P_x$.\par
\end{itemize}
Given observation $\omega$, our goal is to recover $w=Bx$, where $x$ is the signal underlying the observation. A candidate estimate is a Borel function $\widehat{w}(\omega)$ taking value in $\bR^\nu$, and we quantify the performance of such an estimate by its {\sl $(\epsilon,\|\cdot\|)$-risk}
$$
\Risk_{\epsilon,\|\cdot\|}[\widehat{w}|\cX]=\sup_{x\in\cX} \inf\left\{\rho: \Prob_{\xi_x\sim P_x}\{\|Bx-\widehat{w}(Ax+\xi_x)\|>\rho\}\leq\epsilon\;\forall x\in\cX\right\},
$$
that is, the worst, over $x\in\cX$, $(1-\epsilon)$-quantile, taken w.r.t. $P_x$, of the $\|\cdot\|$-magnitude of the recovery error.
\paragraph{Notation.} In the sequel, given a convex compact set, say, $\cY$, in $\bR^n$, we denote by $\cY_\s$ its symmeterization:
$$
\cY_\s=\half(\cY-\cY).
$$
Note that whenever $\cY$ is symmetric w.r.t. the origin, we have $\cY_\s=\cY$. We use ``MATLAB style'' of vector/matrix notation: whenever $H_1,...,H_k$ are matrices of appropriate dimensions, $[H_1,...,H_k]$ stands for horizontal and $[H_1,...,H_k]$ -- for their vertical concatenation. $\bS^n$ stands for the space of $n\times n$ real symmetric matrices equipped with the Frobenius inner product; $\bS^n_+$ is the cone of positive semidefinite matrices from $\bS_n$. Relation $A\succeq B$ ($\Leftrightarrow B\preceq A$) means that $A$ and $B$ are real symmetric matrices of common size such that $A-B$ is positive semidefinite, while $A\succ B$ ($\Leftrightarrow B\prec A$) means that $A$, $B$ are real symmetric matrices of common size such that $A-B$ is positive definite.

\subsection{Generic polyhedral estimate}\label{polyhedest}
A generic polyhedral estimate is as follows:
\begin{quote}
Given the data $A\in \bR^{m\times n},B\in\bR^{\nu\times n},\cX\subset\bR^n$ of the estimation problem stated in Section \ref{theproblem} \ and a ``reliability tolerance'' $\epsilon\in(0,1)$, we somehow specify a positive integer $N$ along with
$N$ linear forms $h_\ell^Tz$ on the observation space $\bR^m$. These forms define linear forms $g_\ell^Tx:=h_\ell^TAx$
on the space of signals $\bR^n$.
Assume that vectors $h_\ell$ are selected in such a way that
\begin{equation}\label{plhdbound}
\forall (x\in\cX): \Prob\{|h_\ell^T\xi_x|>1\}\leq\epsilon/N.
\end{equation}
When setting $H=[h_1,...,h_N]$ (in the sequel, $H$ is referred to as {\sl contrast matrix}),
we clearly have for all $x\in\cX$:
\begin{equation}\label{plhdeq1}
\Prob\left\{\|H^T\omega-H^TAx\|_\infty>1\right\}\leq\epsilon.
\end{equation}
With the polyhedral estimation scheme, we act {\sl as if} all information about $x$ contained in our observation $\omega$ were represented by $H^T\omega$, and we estimate $w=Bx$ by $\wh w=B\wh{x}$, where $\wh{x}=\wh{x}(\omega)$ is a (whatever) vector from $\cX$ compatible with this information, i.e. a solution $\wh{x}$ to the feasibility problem
$$
\hbox{find\ } \;
u\in\cX \;\hbox{such that}\;\|H^T\omega-H^TA{u}\|_\infty\leq 1.
$$
Note that $\wh{x}$ not always is well defined: the above feasibility problem may be unsolvable with positive probability (in fact, with probability $\leq \epsilon$, since, by construction, the true signal $x$ underlying observation $\omega$ is feasible with probability $1-\epsilon$). To circumvent this difficulty, we define $\wh{x}$ according to
\begin{equation}\label{plhdeq2}
\wh{x}^H\in\Argmin_u\left\{\|H^T\omega-H^TAu\|_\infty:u\in\cX\right\}
\end{equation}
so that $\wh{x}^H$ is always well defined and belongs to $\cX$,
and estimate $w$ by $\wh w^H=B\wh{x}^H$.
\end{quote}

\par We have the following immediate observation:
\begin{proposition}\label{polyhedprop} In the situation in question, 
given a contrast matrix $H=[h_1,...,h_N]$ with columns satisfying {\em \rf{plhdbound}}, the quantity
\begin{equation} \label{plhdeq10}
\mR[H]:=\max_{z}\left\{\|Bz\|: \|H^TAz\|_\infty\leq2,z\in2\cX_\s\right\}
\end{equation}
is an upper bound on the $(\epsilon,\|\cdot\|)$-risk of the polyhedral estimate $\widehat{w}^H(\cdot)$:
\begin{equation}\label{plhdeq11}
\Risk_{\epsilon,\|\cdot\|}[\widehat{w}^H|\cX]\leq \mR[H].
\end{equation}
\end{proposition}\noindent
\textbf{Proof} is immediate. Let us fix $x\in\cX$, and let $\cE$ be the set of all realizations of $\xi_x$ such that $\|H^T\xi_x\|_\infty\leq 1$, so that $P_x(\cE)\geq 1-\epsilon$ by (\ref{plhdeq1}). Let us fix a realization $\xi\in\cE$ of the observation noise, let $\omega=Ax+\xi$, and let $\wh{x}=\wh{x}^H(Ax+\xi)$. Then $u=x$ is a feasible solution to the optimization problem (\ref{plhdeq2}) with the value of the objective $\leq1$, implying that the value of this objective at the optimal solution $\wh{x}$ to the problem is $\leq1$ as well, so that $\|H^TA[x-\wh{x}]\|_\infty\leq2$.
Besides this, $z=x-\wh{x}\in2\cX_\s$. We see that $z$ is a feasible solution to (\ref{plhdeq10}), whence
\[\|B[x-\wh{x}]\|=\|Bx-\widehat{w}^H(\omega)\|\leq\mR[H].\] It remains to note that the latter relation holds true whenever $\omega=Ax+\xi$ with $\xi\in\cE$, and for any $x\in\cX$ the $P_x$-probability of the latter inclusion is at least $1-\epsilon$. \qed
\bigskip\par\noindent\textbf{What is ahead.}
In what follows our focus will be on answering the following questions underlying the construction of the polyhedral estimate:
\begin{enumerate}
\item Suppose that, given the data of our estimation problem and a tolerance $\delta\in(0,1)$, we can construct a set $\cH_\delta$ of vectors $h\in\bR^m$ satisfying
the relation
\begin{equation}\label{plhdeq14}
\forall (x\in\cX): \Prob\left\{|h^T\xi_x|>1\right\}\leq\delta.
\end{equation}
With our approach, after the number $N$ of columns in a contrast matrix has been selected, we are free to select these columns  $\cH_{\delta}$, with $\delta=\epsilon/N$, $\epsilon$ being a given reliability tolerance of the estimate we are designing.  Thus, the problem of building sets $\cH_\delta$ satisfying \rf{plhdeq14} arises, the larger $\cH_\delta$, the better.
\item The upper bound $\mR[H]$ on the $(\epsilon,\|\cdot\|)$-risk of the polyhedral estimate $\widehat{w}^H$ is, in general, difficult to compute -- this is the maximum of a convex function over a computationally tractable convex set. Thus, we need to provide computationally efficient upper bounding of
    $\mR[\cdot]$.
\item Finally, given the ``raw materials'' -- set $\cH_\delta$ and an efficiently computable upper bound on the risk of a candidate polyhedral estimate - how to design the best, in terms of (the upper bound on) its risk, polyhedral estimate?
\end{enumerate}
We are about to consider these questions one by one.
\subsection{Specifying sets $\cH_\delta$ for basic observation schemes}\label{setscH}
To specify sets $\cH_\delta$ we are to make assumptions on the distributions of observation noise we want to handle. For the sake of conciseness, in the sequel we restrict ourselves with 3 special observation schemes (below called ``cases'') as follows:
\begin{itemize}
\item {\sl Sub-Gaussian case:} For every $x\in\cX$, the observation noise $\xi_x$ is sub-Gaussian with parameters $(0,\sigma^2I_m)$, $\sigma>0$ (denoted $\xi_x\sim \SG(0,\sigma^2 I_m)$).

Let us denote
\[
\pi_G(h)=\vartheta_{G}\|h\|_2\;\;\mbox{where}\;\; \vartheta_{G}=\sigma\sqrt{2\ln(2/\delta)}.
\]
In the sub-Gaussian case we set
\be
\cH_\delta=\cH_\delta^G:=\{h:\;\pi_G(h)\leq1\}.
\ee{eq:gausset}
\item {\sl Discrete case:} $\cX$ is a convex compact subset of the probabilistic simplex $\mDelta_n=\{x\in\bR^n:x\geq0,\sum_ix_i=1\}$, $A$ is column-stochastic matrix, and
    $$
    \omega={1\over K}\sum_{k=1}^K\zeta_k
    $$
with independent across $k\leq K$ random vectors $\zeta_k$, with $\zeta_k$ taking values $e_i$ with probabilities $[Ax]_i$, $i=1,....,m$, $e_i$ being the basic orths in $\bR^m$.

In this case we put
\bse
\pi_D(h)=2\sqrt{\vartheta_D\max_{x\in\cX}\sum_i[Ax]_ih_i^2 +\mbox{$16\over 9$}\vartheta^2_D\|h\|^2_\infty}\;\;\mbox{with}\;\;\vartheta_D={\ln(2/\delta)\over K},
\ese
and
\be
\cH_\delta=\cH_\delta^D:=\{h:\;\pi_D(h)\leq1\}.
\ee{eq:diset}
\item {\sl Poisson case:} $\cX$ is a convex compact subset of the nonnegative orthant $\bR^n_+$, $A$ is entrywise nonnegative, and the observation $\omega$ stemming from $x\in\cX$ is random vector with independent across $i$ entries $\omega_i\sim\hbox{Poisson}([Ax]_i)$.
In the Poisson case we set
\[
\pi_P(h)=2\sqrt{\vartheta_P \max_{x\in\cX}\sum_i[Ax]_ih_i^2 +\mbox{$4\over 9$} \vartheta_P^2\|h\|^2_\infty}\;\;\mbox{with}\;\;\vartheta_P=\ln(2/\delta),
\]
and
\be
\cH_\delta=\cH_\delta^P:=\{h:\pi_P(h)\leq1\}.
\ee{eq:poiset}
\end{itemize}
We verify in Section \ref{PoissProof} that the sets $\cH^G_\delta,\;\cH^D_\delta$ and $\cH^P_\delta$ as given by \rf{eq:gausset}--\rf{eq:poiset} indeed satisfy
\[
\forall (h\in \cH_\delta,\;x\in \cX)\;\;\;\Prob_x\{|h^T\xi_x|\geq 1\}\leq \delta,
\]
provided that the observation noises $\xi_x$, $x\in\cX$, stem from the respective observation schemes.

\section{Efficient upper-bounding of $\mR[H]$ and contrast design, I.}\label{secteffI}
The scheme for upper-bounding $\mR[H]$ to be presented in this section 
is inspired by the motivating example from the introduction. Indeed, there is a special case of (\ref{plhdeq10}) where $\mR[H]$ is easy to compute -- the case when $\|\cdot\|$ is the uniform norm $\|\cdot\|_\infty$,  whence
$$
\mR[H]=\widehat{\mR}[H]:=2\max_{i\leq \nu}\max_x\left\{\Row_i^T[B]x:x\in\cX_\s,\|H^TAx\|_\infty\leq 1\right\}
$$
is just the maximum of $\nu$ efficiently computable convex functions. It turns out that when $\|\cdot\|=\|\cdot\|_\infty$, it is easy not only to compute $\mR[H]$, but to optimize this risk bound in $H$
as well.
These observations underly the forthcoming developments in this section: under appropriate assumptions, we bound the risk of a polyhedral estimate stemming from a contrast matrix $H$ via the efficiently computable quantity $\widehat{\mR}[H]$ and then show that the resulting risk bounds can be efficiently optimized w.r.t. $H$.  We shall also see that in some simple situations which allow for analytical analysis, like the one in the motivating example, the resulting estimates turn out  to be nearly minimax optimal.
\paragraph{Assumptions.}
 We continue to stay within the setup introduced in Section \ref{theproblem} which we now augment with the following assumptions:
\begin{enumerate}
\item[\textbf{A.1.}] $\|\cdot\|=\|\cdot\|_r$ with $r\in[1,\infty]$.
\item[\textbf{A.2.}] We have at our disposal a sequence $\gamma=\{\gamma_i>0,\,1\leq i\leq \nu\}$ and $\rho\in[1,\infty]$ such that the image of $\cX_\s$ under the mapping $x\mapsto Bx$ is contained in the ``scaled $\|\cdot\|_\rho$-ball''
\begin{equation}\label{plhdEqn1}
\cY=\{y\in\bR^\nu:\|\Diag\{\gamma\}y\|_\rho\leq 1\}.
\end{equation}
\end{enumerate}
\subsection{Simple observation}\label{plhddesignrev}
Let $B_\ell^T$ be $\ell$-th row in $B$, $1\leq \ell\leq \nu$. Let us make the following observation:
\begin{proposition}\label{plhdprop1a}
In the situation described in Section \ref{theproblem}, assuming that Assumptions \textbf{A.1-2} hold, let $\epsilon\in(0,1)$ and a positive real $N\geq\nu$ be given, and let $\pi(\cdot)$ be a norm on $\bR^m$ such that
\begin{equation}\label{plhd20}
\forall (h:\;\pi(h)\leq1,x\in\cX): \Prob\{|h^T\xi_x|>1\}\leq \epsilon/N.
\end{equation}
Let, next, a matrix $H=[H_1,...,H_\nu]$ with $H_\ell\in\bR^{m\times m_\ell}$, $m_\ell\geq1$, and positive reals $\varsigma_\ell$, $1\leq \ell\leq\nu$, satisfy
the relations
\begin{equation}\label{plhdEq700}
\begin{array}{ll}
(a)&\pi(\Col_j[H])\leq1,\,1\leq j\leq N;\\
(b)&\max_x\left\{B_\ell^Tx:\;x\in\cX_\s,\,\|H_\ell^TAx\|_\infty\leq1\right\}\leq\varsigma_\ell,\,1\leq \ell\leq \nu.\\
\end{array}
\end{equation}
Then the quantity $\mR[H]$ as defined in  {\em (\ref{plhdeq10})} can be upper-bounded as follows:
\begin{equation}\label{plhdEq777}
\begin{array}{c}
\mR[H]\leq \Psi(\varsigma):=2\max_v\left\{\|[v_1/\gamma_1;...;v_\nu/\gamma_\nu]\|_r: \,\|v\|_\rho\leq 1,\,0\leq v_\ell\leq \gamma_\ell\varsigma_\ell,\,1\leq \ell \leq \nu\right\}.
\end{array}
\end{equation}
This combines with Proposition \ref{polyhedprop} to imply that
\begin{equation}\label{plhdEq999}
\Risk_{\epsilon,\|\cdot\|}[\wh{w}^H|\cX]\leq\Psi(\varsigma).
\end{equation}
Function $\Psi$ is  nondecreasing on the nonnegative orthant and is easy to compute.
\end{proposition}
\textbf{Proof.} Let $z=2\bar{z}$ be a feasible solution to (\ref{plhdeq10}), so that $\bar{z}\in \cX_\s$ and $\|H^TA\bar{z}\|_\infty\leq 1$. Let $y=B\bar{z}$, so that $y\in\cY$ (see (\ref{plhdEqn1})) due to $\bar{z}\in\cX_\s$ and \textbf{A.2}. Thus,
$\|\Diag\{\gamma\}y\|_p\leq1$.
Besides this, by (\ref{plhdEq700}.$b$) relations $\bar{z}\in \cX_\s$ and $\|H^TA\bar{z}\|_\infty\leq 1$ combine with the symmetry of $\cX_\s$ to imply that
$$
|y_\ell|=|B_\ell^T\bar{z}|\leq \varsigma_\ell,\,\ell\leq\nu.
$$
Taking into account that $\|\cdot\|=\|\cdot\|_r$ by \textbf{A.1}, we see that
$$
\begin{array}{rcl}
\mR[H]&=&\max_z\left\{\|Bz\|_r:z\in 2\cX_\s,\|H^TAz\|_\infty\leq 2\right\}\\
&\leq&2\max_y\left\{\|y\|_r: |y_\ell|\leq\varsigma_\ell,\ell\leq\nu\ \&\ \|\Diag\{\gamma\}y\|_\rho\leq 1\right\}\\
&=&2\max_{v}\left\{\|[v_1/\gamma_1;...;v_\nu/\gamma_\nu]\|_r: \|v\|_\rho\leq 1,0\leq v_\ell\leq \gamma_\ell\varsigma_\ell,\ell\leq \nu\right\},\\
\end{array}
$$
as stated in (\ref{plhdEq777}).\par
It is evident that $\Psi$ is nondecreasing on the nonnegative orthant. Computation of $\Psi$ can be carried out as follows:
\begin{enumerate}
\item When $r=\infty$, we need to compute $\max_{\ell\leq\nu}\max_v\{v_\ell/\gamma_\ell:\|v\|_\rho\leq1,0\leq v_j\leq\gamma_j\varsigma_j,j\leq \nu\}$, so that evaluating $\Psi$ reduces to solving $\nu$ simple convex optimization problems;
\item When $\rho=\infty$, we clearly have $\Psi(\varsigma)=\|[\bar{v}_1/\gamma_1;...;\bar{v}_\nu/\gamma_\nu]\|_r$, $\bar{v}_\ell=\min[1,\gamma_\ell\varsigma_\ell]$;
\item When $1\leq r,\rho<\infty$, passing from variables $v_\ell$ to variables $u_\ell=v_\ell^\rho$, we get
$$
\Psi^r(\varsigma) = 2^r\max_u\left\{\sum_\ell\gamma_\ell^{-r}u_\ell^{r/\rho}:\;\sum_\ell u_\ell\leq1, \; 0\leq u_\ell\leq (\gamma_\ell\varsigma_\ell)^\rho\right\}.
$$
When $r\leq \rho$, the problem on the right hand side is an easily solvable problem of maximizing a simple concave function over a simple convex compact set. When $\infty>r>\rho$, this problem can be solved by Dynamic Programming. \qed
\end{enumerate}
\subsection{Specifying contrasts}\label{speccon}
Risk bound (\ref{plhdEq999}) allows for a straightforward design of contrast matrices. Recalling that $\Psi$ is monotone on the nonnegative orthant, all we need is to select $h_\ell$'s satisfying (\ref{plhdEq700}) and resulting in the smallest possible $\varsigma_\ell$'s, which is what we are about to do now.
\bigskip\par\noindent\textbf{Preliminaries.} Given a vector $b\in\bR^m$ and a norm $s(\cdot)$ on $\bR^m$, consider convex-concave saddle point problem
$$
\Opt=\inf_{g\in\bR^m}\max_{x\in\cX_\s} \left\{\phi(g,x):=[b-A^Tg]^Tx + s(g)\right\}\eqno{(SP)}
$$
along with the induced primal and dual problems
$$
\begin{array}{rcl}
\Opt(P)&=&\inf_{g\in\bR^m}\left[\overline{\phi}(g):=\max_{x\in\cX_\s} \phi(g,x)\right]\\
&=&\inf_{g\in\bR^m}\left[ s(g)+\max_{x\in\cX_\s}[b-A^Tg]^Tx\right]
\end{array}\eqno{(P)}
$$
and
$$\begin{array}{rcl}
\Opt(D)&=&\max_{x\in\cX_\s}\left[\underline{\phi}(g):=\inf_{g\in\bR^m} \phi(g,x)\right]\\
&=&\max_{x\in\cX_\s}\left[\inf_{g\in\bR^m}\left[b^Tx-[Ax]^Tg+s(g)\right]\right]\\
&=&\max_{x}\left[b^Tx:\;x\in\cX_\s,\; q(Ax)\leq1\right]
\end{array}\eqno{(D)}
$$
where $q(\cdot)$ is the norm conjugate to $s(\cdot)$ (we have used the evident fact that $\inf_{g\in\bR^m} [f^Tg+s(g)]$ is either $-\infty$ or 0 depending on whether $q(f)>1$ or $q(f)\leq1$).
Since $\cX_\s$ is compact, we have $\Opt(P)=\Opt(D)=\Opt$ by the Sion-Kakutani theorem. Besides this, $(D)$ is solvable (this is evident) and $(P)$ is solvable as well, since $\overline{\phi}(g)$
is continuous due to the compactness of $\cX_\s$,
and $\overline{\phi}(g)\geq s(g)$, so that $\overline{\phi}(\cdot)$ has bounded level sets. Let $\bar{g}$ be an optimal solution to $(P)$, let $\bar{x}$ be an optimal solution to $(D)$, and let $\bar{h}=\bar{g}/s(\bar{g})$, so that $s(\bar{h})=1$ and $\bar{g}=s(\bar{g})\bar{h}$.
Now let us make the observation as follows:
\begin{proposition}\label{plhdobstolya} In the situation in question, we have
\begin{equation}\label{plhdEqtolya2}
\max_x\left\{|b^Tx|:x\in\cX_\s,|\bar{h}^TAx|\leq1\right\}\leq\Opt.
\end{equation}
In addition, for any matrix $G=[g^1,...,g^M]\in\bR^{m\times M}$ with $s(g^j)\leq1$, $1\leq j\leq M$, one has
\begin{equation}\label{plhdEqtolya3}
\max_x\left\{|b^Tx|:x\in\cX_\s, \|G^TAx\|_\infty\leq1\right\}=\max_x\left\{b^Tx:x\in\cX_\s, \|G^TAx\|_\infty\leq1\right\}\geq\Opt.
\end{equation}
\end{proposition}
\textbf{Proof.} Let $x$ be a feasible solution to the problem in the left hand side of (\ref{plhdEqtolya2}). Replacing, if necessary, $x$ with $-x$, we can assume that $|b^Tx|=b^Tx$. We now have
$$
\begin{array}{l}
|b^Tx|=b^Tx=[\bar{g}^TAx-s(\bar{g})] +\underbrace{[b-A^T\bar{g}]^Tx+s(\bar{g})}_{\leq\overline{\phi}(\bar{g})=\Opt(P)}\leq\Opt(P)+
[s(\bar{g})\bar{h}^TAx-s(\bar{g})]\\
\leq \Opt(P)+s(\bar{g})\underbrace{|\bar{h}^TAx|}_{\leq1}-s(\bar{g}) \leq \Opt(P)=\Opt,\\
\end{array}
$$
as claimed in (\ref{plhdEqtolya2}). Now, the equality in (\ref{plhdEqtolya3}) is due to the symmetry of $\cX_\s$ w.r.t. the origin. To verify the inequality in (\ref{plhdEqtolya3}),
note that $\bar{x}$ satisfies the relations $\bar{x}\in\cX_\s$ and $q(A\bar{x})\leq1$, implying, due to the fact that the columns of $G$ are of $s(\cdot)$-norm $\leq1$,
that $\bar{x}$ is a feasible solution to optimization problems in (\ref{plhdEqtolya3}). As a result, the second quantity in (\ref{plhdEqtolya3}) is at least $b^T\bar{x}=\Opt(D)=\Opt$,
and (\ref{plhdEqtolya3}) follows. \qed
\bigskip\par\noindent\textbf{Designing contrasts.} Propositions \ref{plhdprop1a} and \ref{plhdobstolya} allow for a straightforward solution of the associated contrast design problem, at least in the case of Sub-Gaussian, Discrete, and Poisson observation schemes. Indeed, in these cases, when designing a contrast matrix with $N$ columns, we are supposed to select its columns in the respective sets $\cH_{\epsilon/N}$, see Section \ref{setscH}. Note that these sets are ``nearly independent'' of $N$, because the norms $\pi_G$, $\pi_D$, $\pi_P$ in the description of the respective sets $\cH^G_\delta$, $\cH^D_\delta$, $\cH^P_\delta$ depend on $1/\delta$ only via logarithmic in $1/\delta$ factors. Thus, we lose nearly nothing when assuming that  $N\geq\nu$. So, let us  act as follows:
\begin{quote}
  {\sl  We set $N=\nu$, specify $\overline{\pi}(\cdot)$ as the norm ($\pi_G$, or $\pi_D$, or $\pi_P$) associated with the observation scheme (Sub-Gaussian, or Discrete, or Poisson) in question and $\delta=\epsilon/\nu$, and solve $\nu$ convex optimization problems
 $$
 \begin{array}{rcl}
 \Opt_\ell&=&\min_{g\in\bR^m}\left[\overline{\phi}_\ell(g):=\max_{x\in\cX_\s} \phi_\ell(g,x)\right]\\
 \phi_\ell(g,x)&=&[B_\ell-A^Tg]^Tx+\overline{\pi}(g)
 \end{array}\eqno{(P_\ell)}
 $$
  Next, we convert optimal solution $g_\ell$ to $(P_\ell)$ into a vector $h_\ell\in\bR^m$  by representing $g_\ell=\overline{\pi}(g_\ell)h_\ell$ with $\overline{\pi}(h_\ell)=1$, and set $H_\ell=h_\ell$.
  As a result, we get an $m\times \nu$ contrast matrix $H=[h_1,...,h_\nu]$ which, taken along with $N=
  \nu$, quantities   \begin{equation}\label{plhdEqtolya10}
  \varsigma_\ell=\Opt_\ell,\,1\leq\ell\leq \nu,
  \end{equation}
  and $\pi(\cdot)\equiv \overline{\pi}(\cdot)$,
  in view of the first claim in Proposition \ref{plhdobstolya} as applied with $s(\cdot)\equiv \overline{\pi}(\cdot)$, satisfies the premise of Proposition \ref{plhdprop1a}.

  Consequently, by Proposition \ref{plhdprop1a} we have
  \begin{equation}\label{plhd22}
  \Risk_{\epsilon,\|\cdot\|}[\wh{w}^H|\cX]\leq \Psi([\Opt_1;...;\Opt_\nu]).
  \end{equation}}
  \end{quote}
  \paragraph{Discussion.} Within the framework set in Proposition \ref{plhdprop1a}, optimality of the outlined contrast design for Sub-Gaussian, Discrete and Poisson observation schemes stems from the second claim in Proposition \ref{plhdobstolya} which  states that when $N\geq\nu$ and the columns of the 
    contrast matrix $H=[H_1,...,H_\nu]$ belong to the set $\cH_{\epsilon/N}$ associated with the observation scheme in question, 
    i.e., the norm $\pi(\cdot)$ in the proposition is the norm $\pi_G$, or $\pi_D$, or $\pi_P$ associated with $\delta=\epsilon/N$,
    the quantities  $\varsigma_\ell$ participating in (\ref{plhdEq700}.$b$) cannot be less than $\Opt_\ell$.
    \begin{quote}{\small Indeed, 
    the norm $\pi(\cdot)$ from Proposition \ref{plhdprop1a} is $\geq$ the norm $\overline{\pi}(\cdot)$ participating in $(P_\ell)$ (since the value of $\epsilon/N$ corresponding to $\pi(\cdot)$ is at most $\epsilon/\nu$), implying,
    by (\ref{plhdEq700}.$a$), that the columns of matrix $H$ obeying the premise of the proposition satisfy the relation  $\overline{\pi}(\Col_j[H])\leq1$. Invoking the second part
    of Proposition \ref{plhdobstolya} with $s(\cdot)\equiv\overline{\pi}(\cdot)$, $b=B_\ell$, and  $G=H_\ell$,  and taking (\ref{plhdEq700}.$b$) into account, we conclude that
    $\varsigma_\ell\geq\Opt_\ell$ for all $\ell$, as claimed.}
    \end{quote} Since the bound on the risk of a polyhedral estimate offered by Proposition \ref{plhdprop1a} is the better the less are $\varsigma_\ell$'s, we see that as far as this bound is concerned,   the outlined design procedure is the best possible, provided $N\geq\nu$.
 \par
  An attractive feature of the contrast  design we have just presented is that  it is completely independent of the entities participating in assumptions \textbf{A.1-2} -- these entities affect theoretical risk bounds of the resulting polyhedral estimate, but not the estimate itself.
  \subsection{Illustration: diagonal case}\label{sec:diagcase}
Let us consider the {\sl diagonal case} of our estimation problem, where
\begin{itemize}
\item $\cX=\{x\in\bR^n: \|Dx\|_\rho\leq1\}$, where $D$ is a diagonal matrix with positive diagonal entries $D_{\ell\ell}=:d_\ell$;
\item $m=\nu=n$, and $A$ and $B$ are diagonal matrices with diagonal entries $0<A_{\ell\ell}=:a_\ell$, $0<B_{\ell\ell}=:b_\ell$;
\item $\|\cdot\|=\|\cdot\|_r$;
\item We are in Sub-Gaussian case, that is, $\xi_x\sim \SG(0,\sigma^2I_n)$ for every $x\in\cX$.
\end{itemize}
\par
Let us implement the approach developed in Sections \ref{plhddesignrev} -- \ref{speccon}.
\begin{enumerate}
\item Given reliability tolerance $\epsilon$, we set
\be
\delta=\epsilon/n,\;\;\vartheta_G:=\sigma\sqrt{2\ln(2/\delta)}=\sigma\sqrt{2\ln(2n/\epsilon)}
\ee{plhd25}
and
\[\cH=\cH_\delta^G=\{h\in\bR^n: \pi_G(h):=\vartheta_G\|h\|_2\leq1\};
\]
\item We solve $\nu=n$ convex optimization problems $(P_\ell)$ associated with $\overline{\pi}(\cdot)\equiv \pi_G(\cdot)$, which is immediate: the resulting contrast matrix is
$$
    H=\vartheta_G^{-1}I_n,
    $$ and
\begin{equation}\label{plhdtolya11}
\Opt_\ell=\varsigma_\ell:=b_\ell\min[\vartheta_G/a_\ell,1/d_\ell].
\end{equation}
\end{enumerate}
\bigskip\par\noindent\textbf{Risk analysis.} The $(\epsilon,\|\cdot\|)$-risk of the resulting polyhedral estimate $\wh{w}(\cdot)$ can be bounded by Proposition \ref{plhdprop1a}. Note that setting
$
\gamma_\ell=d_\ell/b_\ell,\;\;1\leq \ell\leq n,
$ 
we meet assumptions  \textbf{A.1-2}, and the above choice of $H$, $N=n$ and $\varsigma_\ell$ satisfies the premise of Proposition \ref{plhdprop1a}. By this proposition,
\begin{equation}\label{plhdNEq36}
\Risk_{\epsilon,\|\cdot\|_r}[\wh{w}^H|\cX]\leq\Psi:=2\max_v\left\{\|[v_1/\gamma_1;...;v_n/\gamma_n]\|_r: \|v\|_\rho\leq1,0\leq v_\ell \leq \gamma_\ell \varsigma_\ell \right\}.
\end{equation}
Let us work out what happens in the {\sl simple case} where
\begin{equation}\label{plhdNEq37}
\begin{array}{ll}
(a)\quad&1\leq \rho\leq r<\infty\\
(b)&a_\ell /d_\ell \ \hbox{and\ }b_\ell /a_\ell \hbox{\ are nonincreasing in $\ell$}\\
\end{array}
\end{equation}

 \begin{proposition}\label{plhdpropsimp} In the just defined simple case,
 let
${\myn}=n$ when
\[\sum_{\ell=1}^{n}\left(\vartheta_Gd_\ell /a_\ell \right)^\rho\leq1,
\] otherwise let ${\mathfrak{n}}$ be the smallest integer such that
\[
\sum_{\ell=1}^{\myn}\left(\vartheta_Gd_\ell /a_\ell \right)^\rho>1,
\]
with $\vartheta_G$ given by {\em (\ref{plhd25})}.
Then for the contrast matrix $H=\vartheta_G^{-1}I_n$ one has
\[
\Risk_{\epsilon,\|\cdot\|_r}[\wh{w}^H|\cX]\leq\Psi\leq 2\left[{\sum}_{\ell=1}^{\mathfrak{n}}(\vartheta_Gb_\ell /a_\ell )^r\right]^{1/r}.
\]
\end{proposition}
For proof, see Section \ref{proofsimp}
\bigskip\par\noindent\textbf{Application.}
Let us apply the result of Proposition \ref{plhdpropsimp} to the ``standard case'' (cf., e.g., \cite{donoho1994minimax,donoho1998minimax,grasmair2018variational}) where
\[
0<
\sqrt{\ln(2n/\epsilon)}\sigma\leq 1,\;a_\ell =\ell^{-\alpha},\;b_\ell =\ell^{-\beta},\;d_\ell =\ell^{\delta}\\
\]
with
$
\beta\geq\alpha\geq0,\;\delta\geq0$ and $(\beta-\alpha) r <1.
$ 
In this case, for large enough $n$, namely, for
\begin{equation}\label{NEq44}
n\geq c\vartheta_G^{-{1\over \alpha+\delta+1/\rho}}\qquad\qquad[\vartheta_G=\sigma\sqrt{2\ln(2n/\epsilon)}]
\end{equation}
(here and in what follows, the factors denoted $c$ and $C$ depend solely on $\alpha,\beta,\delta,r,\rho$) we get
$$
{\myn}\leq C\vartheta_G^{-{1\over \alpha+\delta+1/\rho}}
$$
resulting in
\begin{equation}\label{NEq42}
\Risk_{\epsilon,\|\cdot\|_r}[\wh{w}|\cX]\leq C\vartheta_G^{{\beta+\delta}+1/\rho-1/r\over {\alpha+\delta}+1/\rho}.
\end{equation}
Note that when setting $x=\aic{\Diag\{d_1^{-1},...,d_n^{-1}\}}{D^{-1}}y$, $\bar{\alpha}=\alpha+\delta$, $\bar{\beta}=\beta+\delta$ and treating  $y$, rather than $x$, as the signal underlying the observation, we obtain an estimation problem which is similar to the original one, in which $\alpha,\beta,\delta$ and $\cX$ are replaced, respectively, with $\bar{\alpha}$, $\bar{\beta}$, $\bar{\delta}=0$, and $\cY=\{y:\|y\|_\rho\leq1\}$, and $A$, $B$ are replaced with $\overline{A}=\Diag\{\ell^{-\bar{\alpha}},\ell\leq n\}$, and $\overline{B}=\Diag\{\ell^{-\bar{\beta}},\ell\leq n\}$.
When $n$ is large enough, namely, $n\geq \sigma^{-{1\over \bar\alpha+1/\rho}}$, $\cY$ contains the ``coordinate box''
$$
\overline{\cY}=\{x:|x_\ell|\leq {\mathfrak{m}}^{-1/\rho},{\mathfrak{m}}/2\leq \ell\leq {\mathfrak{m}}, x_\ell=0\ \hbox{otherwise}\}
$$
of dimension $\geq {\mathfrak{m}}/2$, where
$$
{\mathfrak{m}}\geq c\sigma^{-{1\over \bar{\alpha}+1/\rho}}.
$$
 Observe that for all $y\in \overline\cY$,  $\|\bar{A}y\|_2\leq C\mathfrak{m}^{-\bar{\alpha}}\|y\|_2$, and $\|\bar{B}y\|_r\geq c {\mathfrak{m}}^{-\bar{\beta}}\|y\|_r$.
This observation, when combined with the Fano inequality, 
implies that for $\epsilon\ll1$ (cf. \cite{donoho1990minimax}) the minimax optimal w.r.t. the family of all Borel estimates $(\epsilon,\|\cdot\|_r)$-risk on the signal set $\overline\cX=D^{-1}\overline\cY\subset \cX$ is at least
 $$
c \sigma^{\bar{\beta}+1/\rho-1/r\over \bar{\alpha}+1/\rho}.
$$
In other words, in this situation, the upper bound (\ref{NEq42}) on the risk of the polyhedral estimate is within a logarithmic in $n/\epsilon$ factor from the minimax risk. In particular, without surprise, in the case of $\beta=0$ the polyhedral estimates attain well known optimal rates \cite{donoho1994minimax,donoho1997universal,grasmair2018variational}.

\section{Efficient upper-bounding of $\mR[H]$ and contrast design, II.}\label{secteffII}
\subsection{Outline}\label{sec:out}
In this section we develop an approach to the design of polyhedral estimates which is an alternative to that discussed in  Section \ref{secteffI}.
Our present strategy can be outlined as follows. Let us denote by
\[
\cB_*=\{u\in \bR^\nu:\,\|u\|_*\leq 1\}
\]
the unit ball of the norm $\|\cdot\|_*$ {\sl conjugate} to the norm $\|\cdot\|$ in the formulation of the
 estimation problem in Section \ref{theproblem}. Assume that we have at our disposal a technique for bounding quadratic forms on the set  $\cB_*\times\cX_\s$, so that there is an efficiently computable convex function $\cM(M)$ on $\bS^{\nu+n}$
such that
\begin{equation}\label{plhdeq33}
\cM(M)\geq\max_{[u;z]\in\cB_*\times \cX_\s}[u;z]^TM[u;z]\,\,\forall M\in\bS^{\nu+n}.
\end{equation}
Note that the upper bound $\mR[H]$, as defined in  \rf{plhdeq10}, on the risk  of a candidate polyhedral estimate $\wh{w}^H$ given by (\ref{plhdeq10}) is nothing but
\begin{equation}\label{plhdstar}
\mR[H]=2\max_{[u;z]}\bigg\{[u;z]^T\underbrace{\left[\begin{array}{c|c}&{1\over 2}B\cr\hline
{1\over 2}B^T&\cr\end{array}\right]}_{B_+}[u;z]: \begin{array}{l}u\in\cB_*, z\in \cX_\s,\\
z^TA^Th_\ell h_\ell^T Az\leq 1,\ell\leq N\\
\end{array}\bigg\}.
\end{equation}
 When $\lambda\in\bR^N_+$, the constraints
$z^TA^Th_\ell h_\ell^T Az\leq 1$ in (\ref{plhdstar}) can be aggregated to yield the quadratic constraint
$$
z^TA^T\Theta_\lambda Az\leq \mu_\lambda,\,\,\Theta_\lambda=H\Diag\{\lambda\}H^T,\,\mu_\lambda=\sum_\ell\lambda_\ell.
$$
Observe that for every $\lambda\geq0$ we have
\begin{equation}\label{plhdeq30}
\mR[H]\leq 2\cM\bigg(\underbrace{\left[\begin{array}{c|c}&{1\over 2}B\cr\hline {1\over 2}B^T&-A^T\Theta_\lambda A\cr\end{array}\right]}_{B_+[\Theta_\lambda]}\bigg)+2\mu_\lambda.
\end{equation}
\begin{quote}
{\small Indeed, let $[u;z]$ be a feasible solution to the optimization problem $(\ref{plhdstar})$ specifying $\mR[H]$. Then
$$
[u;z]^TB_+[u;z]=[u;z]^TB_+[\Theta_\lambda][u;z]+z^TA^T\Theta_\lambda Az;
$$
the first term in the right hand side is $\leq \cM(B_+[\Theta_\lambda])$ since $[u;z]\in\cB_*\times\cX_\s$, and the second term in the right hand side, as we have already seen, is $\leq\mu_\lambda$, and (\ref{plhdeq30}) follows.}
\end{quote}
Now assume that we have at our disposal a computationally tractable cone
$$
\bH\subset\bS^N_+\times\bR_+
$$
satisfying the following assumption
\begin{quote}
\textbf{Assumption C.} {\sl Whenever $(\Theta,\mu)\in\bH$, we can efficiently find an $n\times N$ matrix $H=[h_1,...,h_N]$ and a nonnegative vector $\lambda\in\bR^N_+$ such that
\begin{equation}\label{plhdeq31}
\begin{array}{ll}
(a)\quad\quad& \hbox{\em the columns $h_\ell$ of $H$ satisfy (\ref{plhdbound})}\\
(b)&\Theta=H\Diag\{\lambda\}H^T\\
(c)&\sum_i\lambda_i\leq \mu
\end{array}
\end{equation}}
\end{quote}
The following simple observation is crucial to what follows:
\begin{proposition}\label{plhdIImain} Consider the estimation problem posed in Section \ref{theproblem}, and let efficiently computable convex function $\cM$ and computationally tractable closed convex cone $\bH$ satisfy  {\em (\ref{plhdeq33})}
\aic{(where $\cB_*$ is the unit ball of the norm conjugate to the norm $\|\cdot\|$ in which the recovery error is measured)}{} and Assumption \textbf{C}, respectively. Consider the convex optimization problem
\begin{equation}\label{plhdprbII}
\begin{array}{c}
\Opt=\min_{\tau,\Theta,\mu}\left\{2\tau+2\mu:(\Theta,\mu)\in\bH,\,\cM(B_+[\Theta])\leq\tau\right\}\\
\left[B_+[\Theta]=\left[\begin{array}{c|c}&{1\over 2}B\cr\hline {1\over 2}B^T&-A^T\Theta A\cr\end{array}\right]\right]\\
\end{array}
\end{equation}
Given a feasible solution $(\tau,\Theta,\mu)$ to this problem, by \textbf{C} we can efficiently convert it to $(H,\lambda)$ such that $H=[h_1,...,h_N]$ with $h_\ell$ satisfying
{\em (\ref{plhdbound})} and $\lambda\geq0$ with $\sum_\ell\lambda_\ell\leq\mu$. We have
$$
\mR[H]\leq 2\tau+2\mu,
$$
whence the $(\epsilon,\|\cdot\|)$-risk of the polyhedral estimate $\wh{w}^H$ satisfies the bound
\begin{equation}\label{plhdeq36}
\Risk_{\epsilon,\|\cdot\|}[\wh{w}^H|\cX]\leq 2\tau+2\mu.
\end{equation}
As a result, we can construct efficiently a polyhedral estimate with $(\epsilon,\|\cdot\|)$-risk arbitrarily close to $\Opt$ (and with risk exactly $\Opt$ if problem { \em(\ref{plhdprbII})} is solvable).
\end{proposition}
\textbf{Proof} is readily given by the reasoning preceding the proposition.  Indeed, with $\tau,\Theta,\mu,H,\lambda$ as in the premise of the proposition, the columns $h_\ell$ of $H$ satisfy (\ref{plhdbound}) by \textbf{C}, implying, by Proposition \ref{polyhedprop}, that $\Risk_{\epsilon,\|\cdot\|}[\wh{w}^H|\cX]\leq\mR[H]$. Besides this, \textbf{C} says that for these $H,\lambda$ it holds  $\Theta=\Theta_\lambda$ and $\mu_\lambda\leq\mu$. Therefore, (\ref{plhdeq30}) combines with the constraints of (\ref{plhdprbII}) to imply that
$\mR[H]\leq 2\tau+2\mu$, and (\ref{plhdeq36}) follows by Proposition \ref{polyhedprop}. \qed
\par
\aic{With the approach we are developing now, presumably good polyhedral estimates will be given by}{
The approach to the design of polyhedral estimate we develop in this section amounts to {reducing} the construction of the estimate (i.e., construction of the contrast matrix $H$) to finding} (nearly) optimal solutions to (\ref{plhdprbII}).
Implementing the proposed approach requires devising techniques for constructing cones $\bH$ satisfying \textbf{C} along with efficiently computable functions $\cM(\cdot)$ satisfying (\ref{plhdeq33}). These tasks are the subjects of the  sections to follow.
\subsection{Specifying cones $\bH$}\label{plhdspecif}
We specify cones $\bH$ in the case when the number $N$ of columns in the candidate contrast matrices is $m$ and under the following assumption on the \aic{given reliability tolerance $\epsilon$ and}{} observation scheme in question:
\begin{quote}
\textbf{Assumption D}. {\sl \aic{The observation scheme in question is such that for properly built}{There is a} computationally tractable convex compact subset $Z\subset\bR^m_+$ intersecting $\inter \bR^m_+$\aic{,}{ such that} the norm $\pi(\cdot)$
$$
\pi(h)=\sqrt{\max_{z\in Z}\sum_iz_ih_i^2}
$$ induced by $Z$ satisfies the relation
$$
\pi(h)\leq 1\;\Rightarrow\; \Prob\{|h^T\xi_x|>1\}\leq\epsilon/m\;\;\forall x\in\cX.
$$}
\end{quote}
Note that Assumption \textbf{D} is satisfied for Sub-Gaussian, Discrete, and Poisson observation schemes: according to the results of Section \ref{setscH},
\begin{itemize}
\item in the Sub-Gaussian case, it suffices to take
$$
Z=\{2\sigma^2\ln(2m/\epsilon)[1;...;1]\};
$$
\item in the Discrete case, it suffices to take
\[
Z={4\ln(2m/\epsilon)\over K}A\cX+{64\ln^2(2m/\epsilon)\over 9K^2}\mDelta_m
\] where
\[
A\cX=\{Ax:x\in\cX\},\;\;\mDelta_m=\{y\in\bR^m:y\geq0,\sum_iy_i=1\};
\]
\item finally, in the Poisson case, it suffices to take
$$
Z=\aic{6}{4}\ln(2m/\epsilon)A\cX+\aic{36}{\mbox{$16\over 9$}}\ln^2(2m/\epsilon)\mDelta_m,
$$
with \aic{the same}{} $A\cX$ and $\mDelta_m$ as \aic{in Discrete case}{above}.
\end{itemize}
Note that in all these cases $Z$  only ``marginally'' -- logarithmically -- depends on $\epsilon$ and $m$.
\par
Under Assumption \textbf{D}, the cone $\bH$ can be built as follows:
\begin{itemize}
\item When $\cZ$ is a singleton: $\cZ=\{\bar{z}\}$, so that $\pi(\cdot)$ is a scaled Euclidean norm, we set
$$
\bH=\{(\Theta,\mu)\in\bS^m_+\times\bR_+: \mu\geq \sum_i\bar{z}_i\Theta_{ii}\}.
$$
Given $(\Theta,\mu)\in \bH$, the $m\times m$ matrix  $H$ and $\lambda\in\bR^m_+$ are defined as follows. We set $S=\Diag\{\sqrt{\bar{z}_1},...,\sqrt{\bar{z}_m}\}$ and compute the eigenvalue decomposition of the matrix $S\Theta S$:
$$
S\Theta S=U \Diag\{\lambda\} U^T,
$$
where $U$ is orthogonal, and set $H=S^{-1}U$, thus ensuring that $\Theta=H \Diag\{\lambda\}H^T$. Since $\mu\geq\sum_i\bar{z}_i\Theta_{ii}$, we have $\sum_i\lambda_i=\Tr(S\Theta S)\leq\mu$. Finally,
a column $h$ of $H$ is of the form $S^{-1}f$ with a $\|\cdot\|_2$-unit vector $f$, implying that \[\pi(h)=\sqrt{\sum_i\bar{z}_i[S^{-1}f]_i^2}=\sqrt{\sum_if_i^2}=1,\] so that $h$ satisfies (\ref{plhdbound}) by Assumption \textbf{D}.
\end{itemize}
\begin{itemize}
\item  When $Z$ is not a singleton, we set
\begin{equation}\label{plhdeq50}
\begin{array}{rcl}
\phi(r)&=&\max_{z\in Z} z^Tr,\\
\varkappa&=&6\ln(2\sqrt{3}m^2),\\
\bH&=&\{(\Theta,\mu)\in\bS^m_+\times\bR_+: \mu\geq\varkappa\phi(\diag(\Theta))\},\\
\end{array}
\end{equation}
where $\diag(Q)$ is the diagonal of a (square) matrix $Q$.
Note that $\phi(r)>0$ whenever $r\geq0$, $r\neq0$, since $Z$ contains a positive vector.
\par
The justification of this construction and the efficient (randomized) algorithm for converting a pair $(\Theta,\mu)\in\bH$ into $(H,\lambda)$ satisfying, along with $(\Theta,\mu)$,  Assumption \textbf{C} are given by the following
\begin{lemma}\label{lemThetamu} Let the norm $\pi(\cdot)$ satisfies Assumption D.\\ {\em (i)} Whenever $H$ is an $m\times m$ matrix with columns $h_\ell$ satisfying $\pi(h_\ell)\leq1$ and $\lambda\in\bR^m_+$, we have
$$
(\Theta_\lambda=H\Diag\{\lambda\}H^T,\;\mu=\varkappa\sum_i\lambda_i)\in\bH.
$$
{\em (ii)} Given $(\Theta,\mu)\in\bH$ with $\Theta\neq0$, we find decomposition $\Theta=QQ^T$ with $m\times m$ matrix $Q$ and fix an orthogonal $m\times m$ matrix $V$ with magnitudes of entries not exceeding $\sqrt{2/m}$ (e.g., the orthogonal scaling of the matrix of the cosine transform).

When $\mu>0$, let us set $\lambda={\mu\over m}[1;...;1]\in\bR^m$ and consider the random matrix
$$
H_\chi=\sqrt{{m\over\mu}}Q\Diag\{\chi\}V,
$$
where $\chi$ is the $m$-dimensional Rademacher random vector (i.e., the entries in $\chi$ are independent of each other random variables taking values $\pm1$ with probabilities $1/2$). We have
\begin{equation}\label{plhdeq40}
H_\chi\Diag\{\lambda\}H_\chi^T\equiv \Theta,\,\lambda\geq0,\sum_i\lambda_i=\mu.
\end{equation}
 Moreover, probability of the event
\begin{equation}\label{plhdeq41}
\pi(\Col_\ell[H_\chi])\leq 1\;\;\forall \ell\leq m
\end{equation}
is at least 1/2. Thus, generating independent samples of $\chi$ and terminating with $H=H_\chi$  when the latter matrix satisfies {\em (\ref{plhdeq41})}, we with probability 1 terminate with $(H,\lambda)$ satisfying Assumption \textbf{C}. Moreover, the probability for the procedure to terminate in course of the first $M=1,2,...$ steps is at least $1-2^{-M}$.
\par
When $\mu=0$, we have $\Theta=0$ (since $\mu=0$ implies $\phi(\diag(\Theta))=0$, which with $\Theta\succeq0$ is possible only when $\Theta=0$); thus, when $\mu=0$, we set $H=0_{m\times m}$ and $\lambda=0_{m\times 1}$.
\end{lemma}
For proof, see Section \ref{lemThetamuproof}.
\par
Note that the lemma states, essentially, that the cone $\bH$ is a tight, up to a logarithmic in $m$ factor, inner approximation of the set
$$
\left\{(\Theta,\mu): \exists (\lambda\in\bR^m_+,H\in\bR^{m\times m}): \begin{array}{l}\Theta=H\Diag\{\lambda\}H^T,\\
\pi(\Col_\ell[H])\leq 1,\,\ell\leq m,\\
 \mu\geq\sum_\ell\lambda_\ell\\
 \end{array}\right\}
$$

\end{itemize}
\subsection{Specifying functions $\cM$}\label{thrdcomp}
In this section we focus on computationally efficient upper-bounding of  maxima of quadratic forms over symmetric w.r.t. the origin convex compact sets, our goal being to specify efficiently computable convex function $\cM(\cdot)$ satisfying (\ref{plhdeq33}). \aic{What we intend to use to this end, is a kind of semidefinite relaxation.}{}
\bigskip\par\noindent\textbf{Cones compatible with convex sets.} Given a nonempty convex compact set $\cY\subset\bR^N$, we say that a cone $\bY$ is {\sl compatible} with $\cY$, if
\begin{itemize}
\item $\bY$ is a closed convex computationally tractable cone contained in $\bS^N_+\times\bR_+$\\
\item one has
\begin{equation}\label{plhdEq0}
\forall (V,\tau)\in\bY:\; \max_{y\in\cY} y^TVy\leq\tau
\end{equation}
\item $\bY$  contains a pair $(V,\tau)$ with $V\succ 0$.
\item relations $(V,\tau)\in \bY$ and $\tau'\geq\tau$ imply that $(V,\tau')\in\bY$.\footnote{The latter requirement is ``for free'' -- passing from a computationally tractable closed convex cone $\bY\subset\bS^N_+\times\bR_+$ satisfying (\ref{plhdEq0}) to the cone $\bY^+=\{(V,\tau):\exists \bar{\tau}\leq \tau: (V,\bar{\tau})\in\bY\}$, we get a larger than $\bY$ cone compatible with $\cY$. It will be clear from the sequel that in our context, the larger is a cone compatible with $\cY$, the better, so that this extension makes no harm.}
\end{itemize}
We call a cone $\bY$ {\sl sharp}, if $\bY$  is a closed convex cone contained in  $\bS^N_+\times\bR_+$ and such that
the only pair $(V,\tau)\in\bY$ with $\tau=0$ is the pair $(0,0)$, or, equivalently, a sequence $\{(V_i,\tau_i)\in \bY,i\geq1\}$ is bounded if and only if the sequence $\{\tau_i,i\geq1\}$ is bounded.
\par
Note that whenever the linear span of $\cY$ is the entire $\bR^N$, every compatible with $\cY$ cone is sharp.\par
Observe that {\sl if $\cY\subset\bR^N$ is a nonempty convex compact set and $\bY$ is a cone compatible with a shift $\cY-a$ of $\cY$, then $\bY$ is compatible with $\cY_\s$} (a symmetrisation of $\cY$).
\begin{quote}{\small Indeed, when shifting a set $\cY$,  its symmeterization $\half[\cY-\cY]$ remains intact, so that we can assume that $\bY$ is compatible with $\cY$.
Now  let $(V,\tau)\in\bY$ and $y,y'\in \cY$. We have
$$
[y-y']^TV[y-y']+\underbrace{[y+y']^TV[y+y']}_{\geq0} =2[y^TVy+[y']^TVy']\leq 4\tau,
$$
whence for $z=\half[y-y']$ it holds $z^TVz\leq\tau$. Since every $z\in\cY_\s$ is of the form $\half[y-y']$ with $y,y'\in\cY$, the claim follows.
\par Note that the claim can be ``nearly inverted:'' {\sl if $0\in\cY$ and $\bY$ is compatible with $\cY_\s$, then the ``widening'' of $\bY$ -- the  cone
$$
\bY^+=\{(V,\tau):\; (V,\tau/4)\in\bY\}
$$
is compatible with $\cY$}  (evident, since when $0\in\cY$, every vector from $\cY$ is proportional, with coefficient 2, to a vector from $\cY_\s$).}
\end{quote}
\bigskip\par\noindent\textbf{Constructing functions $\cM$.} The role of  compatibility in our context becomes clear from the following observation:
{ \begin{proposition}\label{plhdprop2} In the situation described in Section \ref{theproblem}, assume that we have at our disposal cones $\bX$ and $\bU$ compatible, respectively, with $\cX_\s$ and with the unit ball
$$
\cB_*=\{v\in\bR^\nu:\|u\|_*\leq1\}
$$
of the norm $\|\cdot\|_*$ conjugate to $\|\cdot\|$. Given $M\in\bS^{\nu+n}$, let us set
\begin{equation}\label{plhdcM}
\cM(M)=\inf_{X,t,U,s}\left\{t+s:\;(X,t)\in\bX,(U,s)\in\bU, \,\Diag\{U,X\}\succeq M\right\}
\end{equation}
Then $\cM$ is a real-valued efficiently computable convex function on $\bS^{\nu+n}$  such that {\em (\ref{plhdeq33})} takes place: for every $M\in\bS^{n+\nu}$ it holds
$$
\cM(M)\geq \max_{[u;z]\in\cB_*\times\cX_\s}[u;z]^TM[u;z].
$$
In addition, when $\bX$ and $\bU$ are sharp, the infimum in {\em (\ref{plhdcM})} is attained.
\end{proposition}} \noindent
\textbf{Proof} is immediate. Given that the objective of the optimization problem specifying $\cM(M)$ is nonnegative on the feasible set, the fact that $\cM$ is real-valued
is equivalent to problem's feasibility, and the latter is readily given by that fact that $\bX$ is a cone containing a pair $(X,t)$ with $X\succ0$ and similarly for $\bU$. Convexity of $\cM$ is evident. To verify (\ref{plhdeq33}), let $(X,t,U,s)$ form a feasible solution to the optimization problem in (\ref{plhdcM}). When $[u;z]\in\cB_*\times\cX_\s$ we have
$$
[u;z]^TM[u;z]\leq u^TUu+z^TXz\leq s+t,
$$
where the first inequality is due to the $\succeq$-constraint in (\ref{plhdcM}), and the second is due to the fact that $\bU$ is compatible with $\cB_*$, and $\bX$ is compatible with $\cX_\s$. Since the resulting inequality  holds true for all feasible solutions to the optimization problem in (\ref{plhdcM}), (\ref{plhdeq33}) follows. Finally, when $\bX$ and $\bU$ are sharp, (\ref{plhdcM}) is a feasible conic problem with bounded level sets of the objective and as such is solvable. \qed
\subsection{Putting things together}\label{sec:together}
\aic{Combining Propositions \ref{plhdprop2} and \ref{plhdIImain}, we arrive at the following recipe for designing presumably good polyhedral estimates:}{The following statement summarizes our second approach  to the design of polyhedral estimate.}
{ \begin{proposition}\label{plhdproprecipe} In the situation of Section \ref{theproblem}, assume that we have at our disposal cones $\bX$ and $\bU$ compatible, respectively,  with $\cX_\s$ and with the unit ball $\cB_*$ of the norm conjugate to $\|\cdot\|$. Given reliability tolerance $\epsilon\in(0,1)$, assume that we have at our disposal a positive integer $N$ and a computationally tractable cone $\bH$ satisfying, along with $\epsilon$, Assumption \textbf{C}. Consider (clearly feasible) convex optimization problem
\begin{equation}\label{plhdprbfinal}
\Opt=\min_{\Theta,\mu,X,t,U,s}\left\{f(t,s,\mu):=2(t+s+\mu):
\begin{array}{l}(\Theta,\mu)\in\bH,(X,t)\in\bX,(U,s)\in\bU\\
\left[\begin{array}{c|c}U&\half B\cr\hline \half B^T&A^T\Theta A+X\cr\end{array}\right]\succeq0\\
\end{array}\right\}
\end{equation}
Given a feasible solution $\Theta,\mu,X,t,U,s$ to {\em \rf{plhdprbfinal}}, due to \textbf{C}, we can convert, in a computationally efficient manner, $(\Theta,\mu)$ into
$(H,\lambda)$ such that the columns of the $m\times N$ contrast matrix $H$ satisfy {\em (\ref{plhdbound})}, $\Theta=H\Diag\{\lambda\}H^T$, and $\mu\geq\sum_\ell\lambda_\ell$. In this case, the $(\epsilon,\|\cdot\|)$-risk of the polyhedral estimate $\wh{w}^H$ satisfies the bound
\begin{equation}\label{plhdeq70}
\Risk_{\epsilon,\|\cdot\|}[\wh{w}^H|\cX]\leq f(t,s,\mu).
\end{equation}
In particular, we can build, in a computationally efficient manner, polyhedral estimates with risks arbitrarily close to $\Opt$ (and with risk $\Opt$, provided that
{\em (\ref{plhdprbfinal})} is solvable).
\end{proposition}}
\noindent
\textbf{Proof.} Let $\Theta,\mu,X,t,U,s$ form a feasible solution to (\ref{plhdprbfinal}). By the semidefinite constraint in (\ref{plhdprbfinal})
we have
$$
0\preceq \left[\begin{array}{c|c}U&-\half B\cr\hline -\half B^T&A^T\Theta A+X\cr\end{array}\right]=
\Diag\{U,X\}-\underbrace{\left[\begin{array}{c|c} &\half B\cr\hline
\half B^T&-A^T\Theta A\cr\end{array}\right]}_{=:M}.
$$
Hence,  for the function $\cM$ defined in (\ref{plhdcM}) one has
$$
\cM(M)\leq t+s.
$$
Since, by Proposition \ref{plhdprop2}, $\cM$  satisfies (\ref{plhdeq33}), invoking Proposition \ref{plhdIImain} we arrive at
$$
\mR[H]\leq 2(\mu+\cM(M))\leq f(t,s,\mu).
$$
By Proposition \ref{polyhedprop}, this implies the target relation (\ref{plhdeq70}). \qed

\subsection{Compatibility}\label{plhdcompat}
\aic{What is crucial for the design of presumably good polyhedral estimates via}{The approach to design of polyhedral estimates utilizing} the recipe described in Proposition \ref{plhdproprecipe}\aic{,
is}{relies upon} our ability to equip convex ``sets of interest'' (in our context, these are the symmeterization $\cX_\s$ of the signal set and the unit ball $\cB_*$ of the norm conjugate to the norm $\|\cdot\|$) with \aic{cones compatible with these sets.}{compatible cones.}\footnote{Recall that we already know how to specify the second element of the construction, the cone $\bH$.}
Below, we discuss two principal sources of such cones, namely (a) spectratopes/ellitopes, and (b) absolute norms.
\par
More examples of compatible cones can be constructed using a ``compatibility calculus.''
Namely, let us assume that we are given a collection of convex sets (operands) and apply to them some basic operation, such as taking a finite intersection, or arithmetic sum, direct or inverse linear image, or convex hull of a finite union of the sets.
It turns out that cones compatible with the results of such operations can be easily (in a fully algorithmic fashion) obtained from the cones compatible with the operands; see Appendix \ref{app:compcalc} for principal calculus rules.
\par
In view of Proposition \ref{plhdproprecipe}, the larger are the cones $\bX$ and $\bU$ compatible with $\cX_\s$ and $\cB_*$, the better -- the wider is the optimization domain in (\ref{plhdprbfinal}) and, consequently, the less is (the best) attainable risk bound. Given convex compact set $\cY\in\bR^N$, the ``ideal'' -- the largest  candidate to the role of the cone compatible with $\cY$ would be
$$
\bY^*=\{(V,\tau)\in\bS^N_+\times\bR_+: \;\tau\geq \max_{y\in \cY} y^TVy\}.
$$
However, this cone is typically intractable; therefore we look for ``as large as possible'' {\sl tractable} inner approximations of $\bY^*$.

\subsubsection{Cones compatible with ellitopes/spectratopes}\label{sec:comp_spectr}
 The ellitopes and spectratopes as introduced in \cite{l2estimation,JudNem2018} are sets $\cX\subset\bR^n$ representable as linear images
$$\cX=M\cY$$
of {\sl basic ellitopes/spectratopes}, that is, {\sl bounded} sets $\cY$ representable as
\begin{equation}\label{ellspectr}
\begin{array}{lrcll}
(a)&\cY&=&\{y\in\bR^N: \exists r\in\cR: y^TR_\ell y\leq r_\ell,\ell\leq L\}&\hbox{[basic ellitope]}\\
&\multicolumn{4}{c}{R_\ell\succeq 0,\aic{k\leq K}{\ell\leq L},}\\
(b)&\cY&=&\{y\in\bR^N: \exists r\in\cR: R_\ell^2[y]\preceq r_\ell I_{d_\ell},\ell\leq L\}&\hbox{[basic spectratope]}\\
&\multicolumn{4}{c}{R_\ell[y]=\sum_{i=1}^Ny_iR^{\ell i},\,\,R^{\ell i}\in\bS^{d_\ell},}\\
\end{array}
\end{equation}
where $\cR\subset\bR^L_+$ is a convex compact set intersecting with $\inter \bR^L_+$ and such that $0\leq r'\leq r\in\cR$ implies
that $r'\in\cR$.
\par
An ellitope/spectratope is a convex compact set symmetric w.r.t. the origin; as shown in \cite{l2estimation,JudNem2018}, the families of ellitopes/spectratopes admit a kind of fully algorithmic ``calculus'' which demonstrates that these families are closed w.r.t. nearly all basic operations preserving convexity, compactness, and symmetry w.r.t. the origin, including taking finite intersections, direct products, linear images, inverse linear images under linear embeddings, and arithmetic summation. The ``raw materials'' this calculus can be applied to include: for ellitopes -- $\|\cdot\|_p$-balls with $p\in[2,\infty]$, and for spectratopes -- the unit ball of the spectral norm in the space of matrices. In addition, every ellitope is a spectratope as well.
\par
The importance of ellitopes/spectratopes in our context stems from the fact that it is easy to point out cones compatible with these sets.
Specifically, it is shown in \cite{JudNem2018}\footnote{To make the paper self-contained, we reproduce the derivations in Section \ref{AppSpecEll}} that if $\cX=M\cY$, with $\cY$ given by (\ref{ellspectr}.$b$), is a spectratope, then the set
\be
\bX=\left\{(V,\tau)\in \bS^n_+\times\bR_+:\exists\Lambda=\{\Lambda_\ell\in\bS^{d_\ell}_+,\ell\leq L\}: M^TVM\preceq\sum\limits_\ell\cR_\ell^*[\Lambda_\ell],\;\phi_\cR(\lambda[\Lambda])\leq \tau\right\}
\ee{plhdEq20}
where \[
[\cR^*_\ell[\Lambda_\ell]]_{ij}=\Tr(R^{\ell i}\Lambda_\ell R^{\ell j}),\,\lambda[\Lambda]=[\Tr(\Lambda_1);...;\Tr(\Lambda_L)],\;\mbox{and}\;
\phi_\cR(\lambda)=\max\limits_{r\in\cR} r^T\lambda,
\]
is a closed convex cone which is compatible with $\cY$. Similarly, when $\cX=M\cY$ with $\cY$ given by (\ref{ellspectr}.$a$), is an ellitope,
 the set
\begin{equation}\label{plhdEq20a}
\bX=\left\{(V,\tau)\in\bS^{n}_+\times\bR_+:\exists \lambda\in\bR^L_+: M^TVM\preceq\sum_\ell\lambda_\ell R_\ell,\;\phi_\cR(\lambda)\leq\tau\right\}
\end{equation}
is a closed convex cone compatible with $\cY$. In both cases, $\bX$ is sharp, provided that the image space of $M$ is the entire $\bR^n$.
\par
Note that in both these cases $\bX$ is a reasonably tight inner approximation of
\[\bX^*=\{(V,\tau)\in\bS^n_+\times\bR_+: \;\tau\geq \max_{x\in \cX} x^TVx\},\] whenever $(V,\tau)\in \aic{\bY^*}{\bX^*}$, we have $(V,\theta\tau)\in\aic{\bY}{\bX}$, with a moderate $\theta$
(namely, $\theta=O(1)\ln\left(2\sum_\ell d_\ell\right)$ in the spectratopic, and $\theta=O(1)\ln(2L)$ in the ellitopic case, see \cite[Proposition 2.1]{JudNem2018} and \cite[Proposition 3.3]{l2estimation}, respectively).

\subsubsection{Compatibility via absolute norms}
\textbf{Preliminaries.} Recall that a norm $p(\cdot)$ on $\bR^N$ is called {\sl absolute}, if $p(x)$ is a function of the vector $\abs[x]:=[|x_1|;...;|x_N|]$ of the magnitudes of entries in $x$. It is well known that an absolute norm $p$ is monotone on $\bR^N_+$, so that $\abs[x]\leq\abs[x']$ implies that $p(x)\leq p(x')$, and that the norm
$$
p_*(x)=\max_{y:\,p(y)\leq1} x^Ty
$$
conjugate to $p(\cdot)$ is absolute along with $p$.
\par
Let us say that an absolute norm $r(\cdot)$ {\sl fits}  an absolute norm $p(\cdot)$ on $\bR^N$, if for every vector $x$ with $p(x)\leq1$ the entrywise square $[x]^2=[x_1^2;...;x_N^2]$ of $x$ satisfies $r([x]^2)\leq1$. For example, the largest norm $r(\cdot)$ which fits the absolute norm $p(\cdot)=\|\cdot\|_s$, $s\in[1,\infty]$,
 is
 $$
 r(\cdot)=\left\{\begin{array}{ll}\|\cdot\|_1,&1\leq s\leq 2\\
 \|\cdot\|_{s/2},&s\geq 2\\
 \end{array}\right.
 $$
An immediate observation is that an absolute norm $p(\cdot)$ on $\bR^N$ can be ``lifted'' to a norm on $\bS^N$, specifically, the norm
\begin{equation}\label{plhdEq200}
p^+(Y)=p([p(\Col_1[Y]);...;p(\Col_N[Y])]): \bS^N\to\bR_+,
\end{equation}
where $\Col_j[Y]$ is $j$th column in $Y$. It is immediately seen that when $p$ is an absolute norm, the right hand side in (\ref{plhdEq200}) indeed is a norm on $\bS^N$ satisfying the identity
\begin{equation}\label{plhdEq201}
p^+(xx^T)=p^2(x),\,x\in \bR^N.
\end{equation}
\par\noindent
\textbf{Absolute norms and compatibility.} Our interest in absolute norms is motivated by the following immediate
{ \begin{proposition}\label{plhdobs2} Let $p(\cdot)$ be an absolute norm on $\bR^N$, and $r(\cdot)$ be another absolute norm which fits $p(\cdot)$, both norms being computationally tractable. These norms give rise to the computationally tractable and sharp closed convex cone
{\small\begin{equation}\label{plhdEq202}
\bP=\bP_{{p,r}}=\left\{(V,\tau)\in\bS^N_+\times \bR_+: \exists (W\in\bS^N,w\in\bR^N_+):\begin{array}{l} V\preceq W+\Diag\{w\}\\
{[p^+]}_*(W)+r_*(w)\leq\tau\\
\end{array}\right\},
\end{equation}}\noindent
where $[p^+]_*(\cdot)$ is the norm on $\bS^N$ conjugate to the norm $p^+(\cdot)$, and $r_*(\cdot)$ is the norm on $\bR^N$ conjugate to the norm $r(\cdot)$, and this cone is compatible with the unit ball of the norm $p(\cdot)$ (and thus -- with any convex compact subset of this ball).
\end{proposition}}\noindent
Verification is immediate. The fact that $\bP$ is a computationally tractable and closed convex cone is evident. Now let $(V,\tau)\in\bP$, so that $V\succeq0$ and $V\preceq W+\Diag\{w\}$ with $[p^+]_*(W)+r_*(w)\leq \tau$. For $x$ with $p(x)\leq1$ we have
\bse
x^TVx&\leq& x^T[W+\Diag\{w\}]x=\Tr(W[xx^T])+w^T[x]^2\\
&\leq& p^+(xx^T)[p^+]_*(W)+r([x]^2)r_*(w)
=p^2(x)[p^+]_*(W)+r_*(w)\\
&\leq& [p^+]_*(W)+r_*(w)\leq\tau,
\ese
whence $x^TVx\leq \tau$ for all $x$ with $p(x)\leq1$. \qed
Let us look at the proposed construction in the case where $p(\cdot)=\|\cdot\|_s$, $s\in[1,\infty]$, and let $r(\cdot)=\|\cdot\|_{\bar{s}}$, $\bar{s}=\max[s/2,1]$. Setting $s_*={s\over s-1}$, $\bar{s}_*={\bar{s}\over\bar{s}-1}$, we clearly have
\[
[p^+]_*(W)=\|W\|_{s_*}:=\left\{\begin{array}{ll} \left(\sum_{i,j}|W_{ij}|^{s_*}\right)^{1/s_*},&s_*<\infty\\
\max_{i,j}|W_{ij}|,&s_*=\infty\\
\end{array}\right.,\,\,r_*(w)=\|w\|_{\bar{s}_*},
\]
resulting in
\begin{equation}\label{plhdEq204}
\begin{array}{rcl}
\bP^s&:=&\bP_{\|\cdot\|_s,\|\cdot\|_{\bar{s}}}\\
&=&\left\{(V,\tau):V\in\bS^N_+,\exists (W\in\bS^N,w\in\bR^N_+):
\begin{array}{l}V\preceq W+\Diag\{w\},\\
\|W\|_{s_*}+\|w\|_{\bar{s}_*}\leq\tau
\end{array}\right\}.
\end{array}
\end{equation}
By Proposition \ref{plhdobs2}, $\bP^s$ is compatible with the unit ball of $\|\cdot\|_s$-norm on $\bR^N$ (and therefore with every closed convex subset of this ball).
\par
When $s=1$, that is, $s_*=\bar{s}_*=\infty$,  (\ref{plhdEq204}) results in
\begin{equation}\label{plhdEq204a}
\begin{array}{rcl}
\bP^1&=&\left\{(V,\tau):V\succeq0,\exists (W\in\bS^N,w\in\bR^N_+):\begin{array}{l} V\preceq W+\Diag\{w\},\\
\|W\|_\infty+\|w\|_\infty\leq\tau\\
\end{array}\right\}\\
&=&\{(V,\tau):V\succeq0,\|V\|_\infty\leq\tau\},\\
\end{array}
\end{equation}
and it is easily seen that the situation is a good as it could be, namely, $$\bP^1=\left\{(V,\tau):V\succeq0,\max_{\|x\|_1\leq 1}x^TVx\leq\tau\right\}.$$
It can be shown (see Section \ref{plhdeasy}) that when $s\in[2,\infty]$, so that $\bar{s}_*={s\over s-2}$, (\ref{plhdEq204}) results in
\begin{equation}\label{plhdEq204b}
\bP^s=\left\{(V,\tau):V\succeq0,\exists (w\in\bR^N_+):\;V\preceq\Diag\{w\}\ \&\  \|w\|_{{s\over s-2}}\leq\tau\right\}.
\end{equation}
Note that
$$
\bP^2=\{(V,\tau):\;V\succeq0,\,\|V\|_{\mathrm{\tiny sp}}\leq\tau\}
$$
where $\|\cdot\|_{\mathrm{\tiny sp}}$ is the spectral norm, and this is {\sl exactly} the largest cone compatible with the unit Euclidean ball.
\par
When $s\geq 2$, the unit ball $\cY$ of the norm $\|\cdot\|_s$ is an ellitope:
{\small$$
\{y\in\bR^N:\|y\|_s\leq 1\}=\{y\in\bR^N:\exists (t\geq0,\|t\|_{\bar{s}}\leq1): y^TR_\ell y:=y_\ell^2\leq t_\ell,\,\ell\leq L=N\},
$$}\noindent
so that \aic{one of the cones}{a} compatible with $\cY$ \aic{}{cone} is given by (\ref{plhdEq20a}) with \aic{the identity matrix in the role of $M$}{$M=I_N$}. It comes as no surprise that, as it is immediately seen, the latter cone is nothing but the cone given by (\ref{plhdEq204b}).

\subsection{Near-optimality of polyhedral estimate in the spectratopic sub-Gaussian case}
As an instructive application of the design of polyhedral estimate  developed in this section, consider the special case of the estimation problem stated in Section \ref{theproblem}, where
\begin{enumerate}
\item The signal set $\cX$ and the unit ball $\cB_*$ of the norm conjugate to $\|\cdot\|$ are spectratopes:
$$
\begin{array}{rcl}
\cX&=&\{x\in\bR^n: \exists  t\in\cT: R_k^2[x]\preceq t_k I_{d_k},\,1\leq k\leq K\},\\
\cB_*&=&\{z\in\bR^\nu: \exists  y\in\cY:z=My\},\\
&&\cY:=\{y\in\bR^q:\exists r\in\cR: S_\ell^2[y]\preceq r_\ell I_{f_\ell},\,1\leq\ell\leq L\}, \\
\end{array}
$$
\item For every $x\in\cX$, observation noise is sub-Gaussian, i.e. $\xi_x\sim \SG(0,\sigma^2I_m)$.
\end{enumerate}
We are about to show that in the present situation, {\sl the polyhedral estimate yielded by the approach described in Sections \ref{sec:out}--\ref{sec:together}, i.e., yielded by the efficiently computable (high accuracy near-) optimal solution to the optimization problem {\em \rf{plhdprbfinal}} is near-optimal in the minimax sense.}

Given reliability tolerance $\epsilon\in(0,1)$,
the recipe for constructing a $m\times m$ contrast matrix $H$ as prescribed by
Proposition \ref{plhdproprecipe} is as follows:
\begin{itemize}
\item {We set}
$$
Z=\{\vartheta^2[1;...;1]\},\,\vartheta=\sigma\kappa,\,\kappa=\sqrt{2\ln(2m/\epsilon)},
$$
and utilize the construction from Section \ref{plhdspecif},
thus arriving at the cone
$$
\bH=\{(\Theta,\mu)\in\bS^m_+\times\bR_+: \sigma^2\kappa^2 \Tr(\Theta)\leq\mu\}
$$
satisfying the requirements of Assumption \textbf{C}.
\item We specify the cones $\bX$ and $\bU$ compatible with $\cX_\s=\cX$, and $\cB_*$, respectively, according to (\ref{plhdEq20}).
\end{itemize}
The resulting problem (\ref{plhdprbfinal}), after immediate straightforward simplifications, reads
\begin{equation}\label{plhdreadsini}
\begin{array}{rcl}
\Opt&=&\min\limits_{\Theta,U,\Lambda,\Upsilon}\bigg\{2\left[\phi_\cR(\lambda[\Upsilon])+\phi_{\cT}(\lambda[\Lambda])+\sigma^2\kappa^2\Tr(\Theta)\right]:\\
&&\begin{array}{l}\Theta\succeq0,\,U\succeq0,\Lambda=\{\Lambda_k\succeq0,k\leq K\},\,\Upsilon=\{\Upsilon_\ell\succeq0,\ell\leq L\},\\
\left[\begin{array}{c|c}U&{1\over 2}B\cr\hline {1\over 2}B^T&A^T\Theta A+\sum_k \cR_k^*[\Lambda_k]\cr\end{array}\right]\succeq0,\;
M^TUM\preceq \sum_\ell\cS^*_\ell[\Upsilon_\ell]\\
\end{array}\Bigg\}
\end{array}
\end{equation}
where
\bse
[\cR_k^*[\Lambda_k]]_{ij}&=&\Tr(R^{k i}\Lambda_kR^{kj}) \quad\left[R_k[x]=\sum_ix_iR^{ki}\right],\\
{[\cS_\ell^*[\Upsilon_\ell] ]}_{ij}
&=&\Tr(S^{\ell i}\Upsilon_\ell S^{\ell j})\quad
\left[S_\ell[u]=\sum_iu_iS^{\ell i}\right],
\ese
and
\[\lambda[\Lambda]=[\Tr(\Lambda_1);...;\Tr(\Lambda_K)],\,\lambda[\Upsilon]=[\Tr(\Upsilon_1);...;\Tr(\Upsilon_L)],\,\phi_W(f)=\max_{w\in W} w^Tf.\]
Let now
\[
\RiskOpt_\epsilon=\inf\limits_{\wh{w} (\cdot)}\sup\limits_{x\in \cX}\inf\left\{\rho:
\Prob_{\xi\sim\cN(0,\sigma^2I)}\{\|Bx-\wh{w}(Ax+\xi)\|>\rho\}\leq \epsilon\,\forall x\in\cX\right\},
\]
be the minimax optimal $(\epsilon,\|\cdot\|)$-risk of estimating $w=Bx$ in Gaussian observation scheme where $\xi_x\sim \cN(0,\sigma^2 I_m)$ independently of $x\in X$.
\begin{proposition}\label{pr:near}
When $\epsilon\leq1/8$, the polyhedral estimate $\wh w^H$ yielded by a feasible near optimal, in terms of the objective, solution to problem {\em (\ref{plhdreadsini})} is minimax optimal within the logarithmic factor, namely
\[
\begin{array}{rcl}
\Risk_{\epsilon,\|\cdot\|}[\widehat{w}|\cX]&\leq&  O(1)\sqrt{\ln\Big(\sum_kd_k\Big)\ln \Big(\sum_\ell f_\ell\Big)\ln(2m/\epsilon)}\; \RiskOpt_{{1\over 8}}\\
&\leq& O(1)\sqrt{\ln\Big(\sum_kd_k\Big)\ln \Big(\sum_\ell f_\ell\Big)\ln(2m/\epsilon)}\; \RiskOpt_{\epsilon},\\
\end{array}
\]
where $O(1)$ is an absolute constant.
\end{proposition}
For proof, see Section \ref{proofnear}.
{
\paragraph{Discussion.}
It is worth mentioning that the approach described in Section \ref{secteffI} is complementary to the approach discussed in this section. In fact, it is easily seen that the bound $\Opt$ for the risk of the polyhedral estimate stemming from \rf{plhdprbfinal} is suboptimal in the simple situation described in the motivating example from Introduction. Indeed, let $\cX$ be the unit $\|\cdot\|_1$-ball, let $\|\cdot\|=\|\cdot\|_2$, and let us consider the problem of estimating $x\in \cX$ from the direct observation $\omega=x+\xi$ with Gaussian observation noise $\xi\sim \cN(0,\sigma^2I)$. We equip the ball $\cB_*=\{u\in \bR^n:\;\|u_2\|_2\leq 1\}$ with the cone
\[\bU=\bP^2=\{(U,\tau):\;U\succeq 0,\,\|U\|_{\mathrm{\tiny sp}}\leq \tau\}
\]
and $\cX$ -- with the cone
\[
\bX=\bP^1=\{(X,t):\;X\succeq 0,\,\|X\|_\infty\leq t\},
\]
(note that both cones are the largest w.r.t. inclusion cones compatible with the respective sets). The corresponding problem \rf{plhdprbfinal} reads
\be
\Opt&=&\min_{\Theta,X,U}\left\{2\Big(\kappa^2\sigma^2\Tr(\Theta)+\max_i X_{ii}+\|U\|_{\mathrm{\tiny sp}}\Big):
\begin{array}{l}\Theta\succeq 0,\;X\succeq 0,\;U\succeq 0,\\
\left[\begin{array}{c|c}U&{\half}I_n\cr\hline {\half}I_n&\Theta+X\cr\end{array}\right]\succeq0\\
\end{array}\right\}\nn
&=&\min_{\Theta,X,U}\left\{2\Big(\kappa^2\sigma^2\Tr(\Theta)+\max_i X_{ii}+\tau\Big):
\begin{array}{l}\Theta\succeq 0,\;X\succeq 0,\;U\succeq 0, \\
\left[\begin{array}{c|c}\tau I_n&{\half}I_n\cr\hline {\half}I_n&\Theta+X\cr\end{array}\right]\succeq0\\
\end{array}\right\}
\ee{NEq10}
Observe that every $n\times n$ matrix of the form $Q=EP$, where $E$ is diagonal with diagonal entries $\pm1$, and $P$ is a permutation matrix, induces a symmetry
$(\Theta,X,\tau)\mapsto (Q\Theta Q^T,QXQ^T,\tau)$ of the second optimization problem in (\ref{NEq10}), that is, a transformation which maps the feasible set onto itself and keeps the objective intact. Since the problem is convex and solvable, we conclude that it has an optimal solution which remains intact under the symmetries in question, i.e., solution with scalar matrices $\Theta=\theta I_n$ and $X=u I_n$. As a result,
\be
\Opt=\min_{\theta\geq0,u\geq0,\tau}\left\{2(\kappa^2\sigma^2n\theta+u+\tau): \;\tau(\theta+u)\geq \four\right\}=
2\min\left[\kappa\sigma\sqrt{n},1\right],
\ee{neq1001}
A similar derivation shows that the value $\Opt$ remains intact if we replace the set $\cX=\{x:\|x\|_1\leq1\}$ with $\cX=\{x:\|x\|_s\leq1\}$, $s\in [1, 2]$, and the cone  $\bX=\bP^1$ with $\bX=\bP^s$, see \rf{plhdEq204}.  Since the $\Theta$-component of an optimal solution to (\ref{NEq10}) can be selected to be scalar,
the contrast matrix  $H$ we end up with can be selected to be the unit matrix. An unpleasant observation is that when $s<2$, the quantity $\Opt$ given by \rf{neq1001} ``heavily overestimates'' the actual 
risk of the polyhedral estimate with $H=I_n$. Indeed, the analysis of this estimate in Section \ref{sec:diagcase} results in the risk bound (up to a logarithmic in $n$ factor) $\min[\sigma^{1-s/2},\sigma\sqrt{n}]$, what can be much less than $\Opt=2\min\left[\kappa\sigma\sqrt{n},1\right]$, e.g., in the case of large $n$, and $\sigma\sqrt{n}=O(1)$.
\par
The good news is that whenever the approaches developed in this section and in Section \ref{secteffI} are applicable, they can be utilized simultaneously.
The underlying observation is that
\begin{quote}
(!) {\em \aic{With our approach to ensuring {\em (\ref{plhdbound})} and  (at least) for the
observation schemes we are considering,} In the problem setting described in Section \ref{sect1}, a collection of $K$ candidate polyhedral estimates  can be assembled into a single polyhedral estimate with the (upper bound on the) risk, as given by Proposition \ref{polyhedprop}, being nearly the minimum of the risks of estimates we aggregate.}
\end{quote}
Indeed, given an observation
scheme (that is, collection of probability distributions $P_x$ of noises $\xi_x$, $x\in\cX$), assume we have at our disposal norms
$$
\pi_\delta(\cdot):\bR^m\to\bR
$$
parameterized by $\delta\in(0,1)$ such that $\pi_\delta(h)$, for every $h$, is the larger the less is $\delta$, and (cf. Section \ref{plhdspecif})
$$
 \forall (x\in \cX,\delta\in(0,1),h\in\bR^m): \pi_\delta(h)\leq 1\Rightarrow \Prob_{\xi\sim P_x}\{\xi:|h^T\xi|>1\}\leq \delta.
$$
 Assume also (as indeed is the case in all our constructions) that we ensure (\ref{plhdbound}) by imposing on the columns  $h_\ell$ of an $m\times N$ contrast matrix $H$ in question the restrictions $\pi_{\epsilon/N}(h_\ell)\leq1$.
\par
Now suppose that given the risk tolerance $\epsilon\in(0,1)$, we have generated somehow $K$ candidate contrast matrices $H_k\in\bR^{m\times N_k}$ such that
$$
\pi_{\epsilon/N_k}(\Col_j[H_k])\leq1,\,j\leq N_k,
$$
so that the $(\epsilon,\|\cdot\|)$-risk of the polyhedral estimate yielded by the contrast matrix $H_k$ does not exceed
$$
\mR_k=\max_x\left\{\|Bx\|:x\in2\cX_\s,\|H_k^TAx\|_\infty\leq 2\right\}.
$$
Let us combine the contrast matrices $H_1,...,H_K$ into a single contrast matrix $H$ with $N=N_1+...+N_K$ columns by normalizing the columns of the matrix $[H_1,...,H_K]$ to have $\pi_{\epsilon/N}$-norms equal to 1, so that
$$
H=[\bar{H}_1,...,\bar{H}_K] ,\,\Col_j[\bar{H}_k]=\theta_{jk}\Col_j[H_k]\,\,\forall (k\leq K,j\leq N_k)
$$
with
$$
\theta_{jk}={1\over \pi_{\epsilon/N}(\Col_j[H_k])}\geq \vartheta_k:=\min_{h\neq0}{ \pi_{\epsilon/N_k}(h)\over \pi_{\epsilon/N}(h)},
$$
where the concluding $\geq$ is due to $ \pi_{\epsilon/N_k}(\Col_j[H_k])\leq1$.
We claim that in terms of $(\epsilon,\|\cdot\|)$-risk, the polyhedral estimate yielded by $H$ is ``\aic{basically}{almost}} as good'' as the best of the polyhedral estimates yielded by the contrast matrices $H_1,...,H_K$, specifically,
$$
\mR[H]:=\max_x\left\{\|Bx\|:x\in 2\cX_\s, \|H^TAx\|_\infty\leq2\right\}\leq \min_k\vartheta_k^{-1}\mR_k.\footnotemark
$$
\footnotetext{This is the precise ``quantitative expression'' of the observation (!).}
The justification is readily given by the following observation: when $\vartheta\in(0,1)$, we have
$$
\mR_{k,\vartheta}:=\max_x\left\{\|Bx\|:x\in2\cX_\s,\|H_k^TAx\|_\infty\leq 2/\vartheta\right\}\leq\mR_k/\vartheta.
$$
Indeed, when $x$ is a feasible solution to the maximization problem specifying $\mR_{k,\vartheta}$, $\vartheta x$ is a feasible solution to the problem specifying $\mR_k$, implying that $\vartheta\|Bx\|\leq\mR_k$.  It remains to note that we clearly have $\mR[H]\leq\min_k\mR_{k,\vartheta_k}$.
\par
The bottom line is that the just described \aic{assembling}{aggregation of} contrast matrices $H_1,...,H_K$ into a single contrast matrix $H$ results in polyhedral estimate which in terms of upper bound $\mR[\cdot]$ on  its $(\epsilon,\|\cdot\|)$-risk is, up to the factor $\bar\vartheta=\max_k\vartheta_k^{-1}$, not worse than the best of the $K$ estimates yielded by the original contrast matrices. Consequently, whenever $\pi_\delta(\cdot)$ grows slowly as $\delta$ decreases, the ``price'' $\bar\vartheta$ of assembling the original estimates is quite moderate. For example, in the case of basic observation schemes  (Sub-Gaussian, Discrete, and Poisson) described in Section \ref{setscH}, $\bar\vartheta$ is logarithmic in $\max_kN_k^{-1}(N_1+...+N_K)$, and  $\bar\vartheta=1+o(1)$ as $\epsilon\to+0$ for $N_1,...,N_K$ fixed.
\subsection{Numerical illustration}
To illustrate the performance of the polyhedral estimates we compare it numerically with a ``presumably good'' linear estimate. Our setup is deliberately
simple: the signal set $\cX$ is just the unit box $\{x\in\bR^n:\|x\|_\infty\leq1\}$, $B\in\bR^{n\times n}$ is ``numerical double integration:'' for a $\delta>0$,
$$
B_{ij}=\left\{\begin{array}{ll} \delta^2(i-j+1),&j\leq i\\
0,&j>i\\
\end{array}\right.,
$$
so that $x$, modulo boundary effects, is the second order finite difference derivative of $w=Bx$:
$$
x_i={w_i-2w_{i-1}+w_{i-2}\over\delta^2},\,2<i\leq n,
$$
and $Ax$ is comprised of $m$ randomly selected entries of $Bx$. The observation is
$$
\omega=Ax+\xi,\,\xi\sim\cN(0,\sigma^2I_m).
$$
and the recovery norm is $\|\cdot\|_2$. In other words, we want to recover a restriction of twice differentiable function of one variable on the $n$-point regular grid on the segment $\Delta=[0,n\delta]$ from noisy observations of this restriction taken along $m$ randomly selected points of the grid.
A priori information on the function is
that the magnitude of its second order derivative does not exceed 1.
\par
Note that in the considered situation both linear estimate  $\wh{w}_H$ yielded by Proposition 3.3 in
\cite{JudNem2018} and polyhedral estimate $\wh{w}^H$  yielded by Proposition \ref{plhdprop2}, are near-optimal in the minimax sense in terms of their $\|\cdot\|_2$-, resp., $(\epsilon,\|\cdot\|_2)$-risk. \aic{The goal of the presented experiment is to compare numerically performance of the estimates on simulated data.}{}
\par
\begin{figure}[ht]
\begin{center}
\begin{tabular}{cc}
\includegraphics[width=0.45\textwidth]{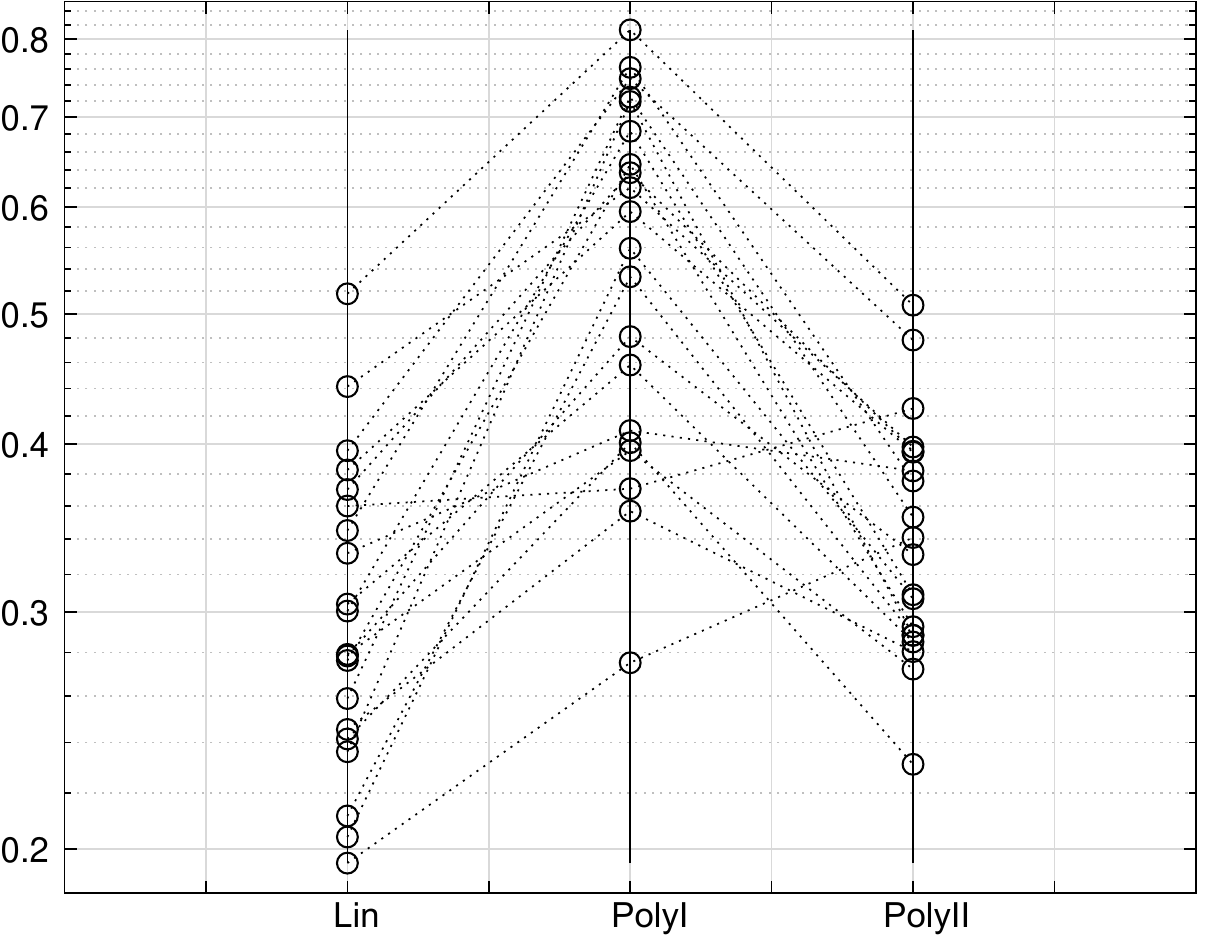}&
\includegraphics[width=0.45\textwidth]{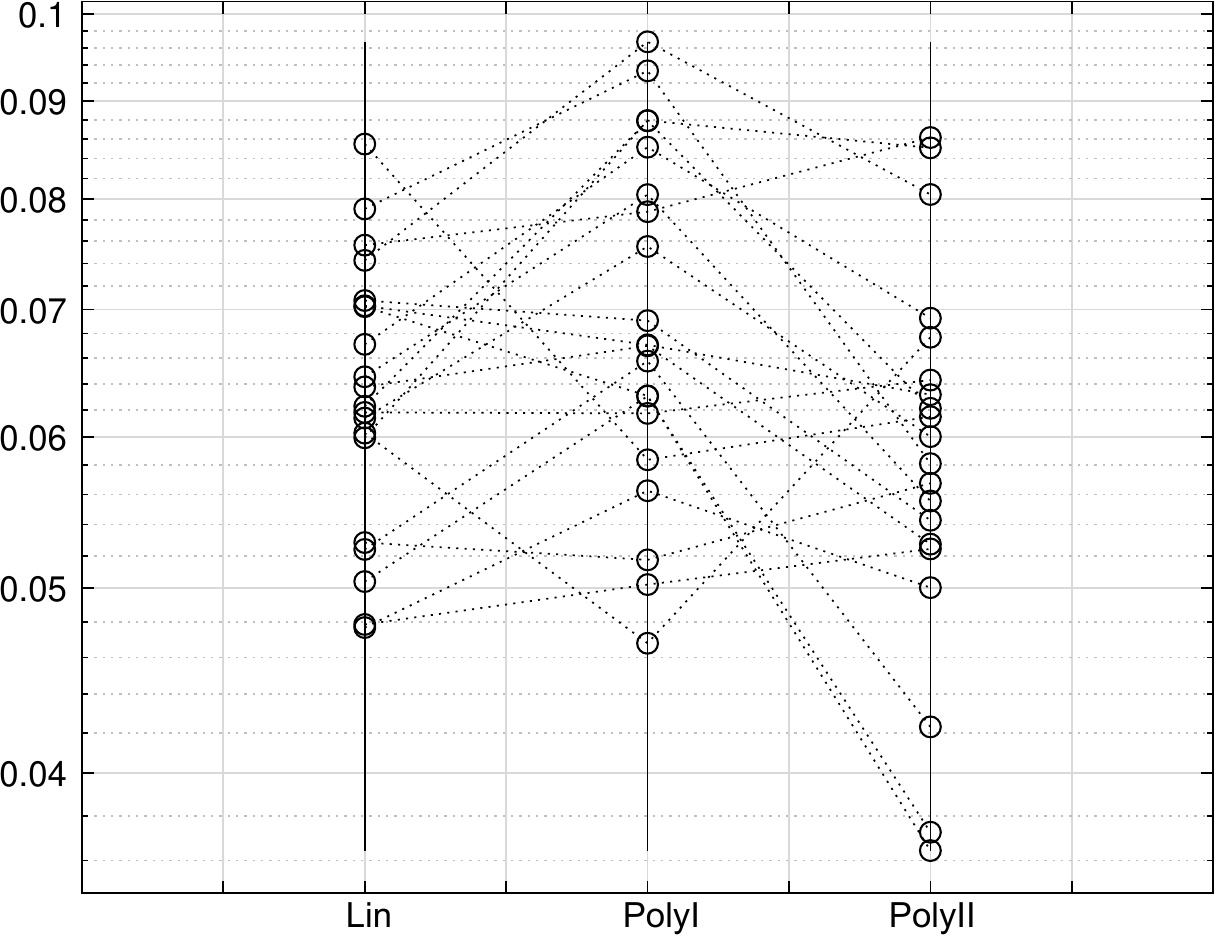}\\
$\sigma=0.1$&$\sigma=0.01$\\
\includegraphics[width=0.45\textwidth]{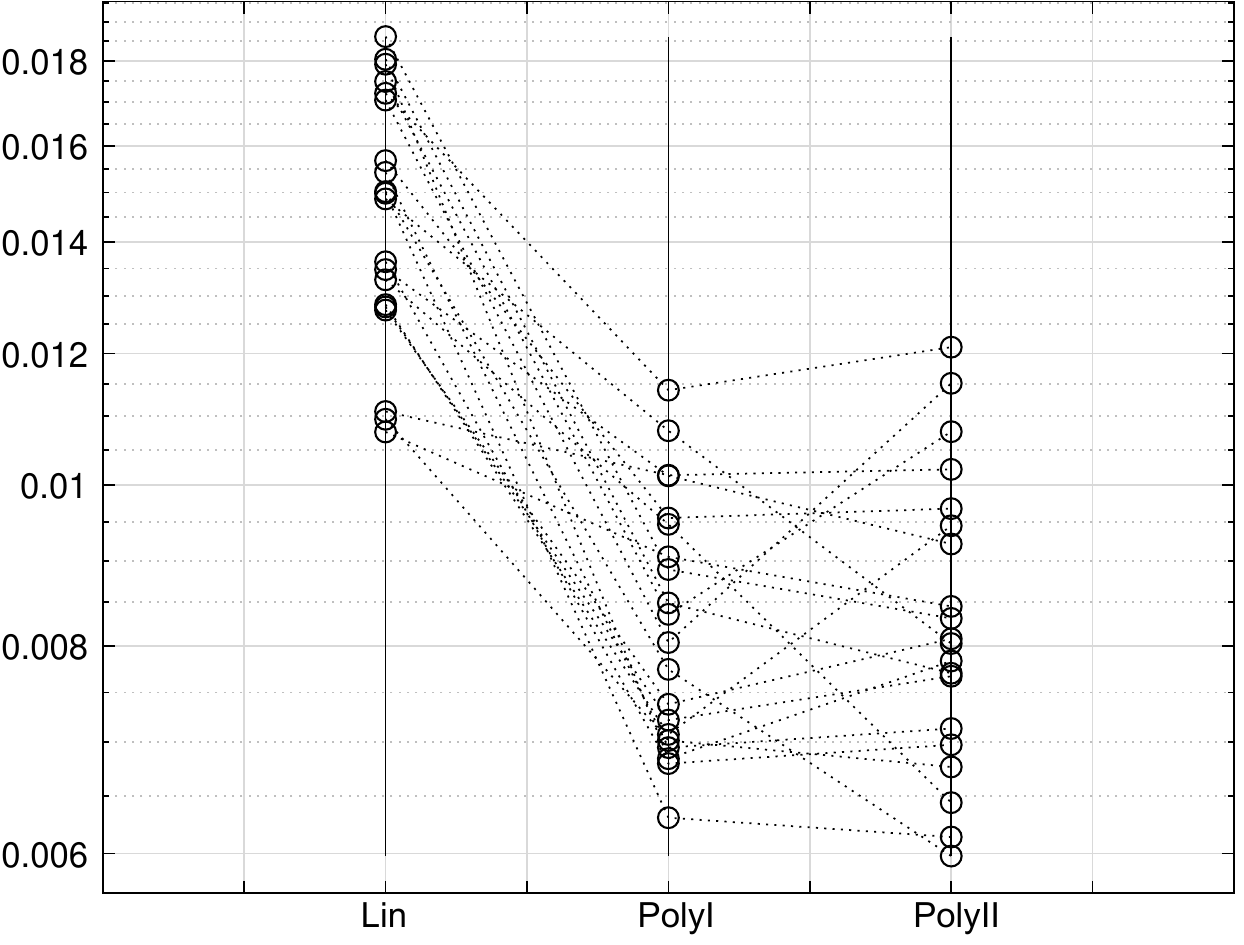}&
\includegraphics[width=0.45\textwidth]{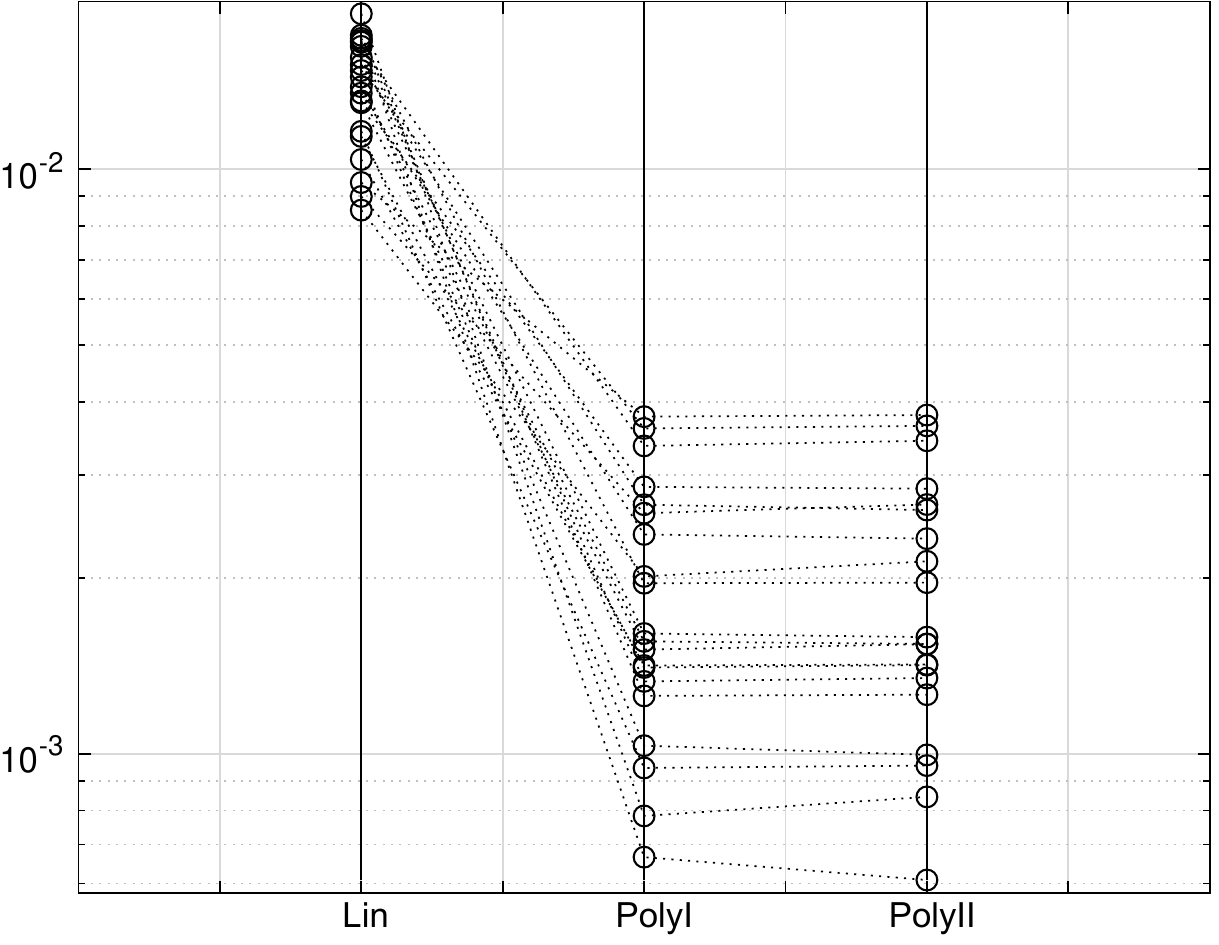}\\
$\sigma=0.001$&$\sigma=0.0001$
\end{tabular}
\caption{\label{plvslnfig}
Recovery errors for near-optimal linear estimate (left) and for polyhedral estimates yielded by Proposition \ref{plhdprop2} ({\em PolyI}, middle) and by the construction from Section \ref{secteffI} ({\em PolyII}, right), 20 simulations per each value of $\sigma$.}
\end{center}
\end{figure}

In the experiments reported in Figure \ref{plvslnfig}, we used $n=64$, $m=32$, and $\delta=4/n$ (i.e., $\Delta=[0,4]$); the reliability parameter for the polyhedral estimate was set to $\epsilon=0.1$. For different noise levels $\sigma=\{0.1,0.01,0.001,0.0001\}$ we generate at random 20 signals $x$ from $\cX$ and record the $\|\cdot\|_2$-recovery errors of the linear and the polyhedral estimates. In addition to testing the nearly optimal polyhedral estimate {\em PolyI} yielded by Proposition \ref{plhdprop2}, we also record the  performance of the polyhedral estimate {\em PolyII} yielded the construction from Section \ref{secteffI}. The observed $\|\cdot\|_2$-recovery errors of the three estimates, sorted in the non-descending order, are plotted in Figure \ref{plvslnfig}.
\par
In these simulations, the three estimates exhibit similar empirical performance. However, when the noise level becomes small, polyhedral estimates seem to outperform the linear one. In addition, the estimate {\em PolyII} seems to ``work'' better than or, at the very worst, similarly to {\em PolyI} in spite of the fact that in the situation in question the estimate {\em PolyI,} in contrast to {\em PolyII,} is provably near-optimal.

\appendix
\section{Appendix}
\subsection{Executive summary on computationally tractable convex sets}\label{comptract}
Detailed description of what ``computationally tractable convex set/function'' is goes beyond the scope of this paper; an interested reader is referred to \cite{LMCO}. Here we restrict ourselves with  ``executive summary'' as follows: ``for all practical purposes'' (including all applications considered in this paper), a computationally tractable convex set is a set $X\subset \bR^N$ given by {\sl semidefinite representation}
$$
X=\{x\in\bR^N:\exists u\in\bR^M: \cA(x,u)\succeq 0\},
$$
where $\cA(\cdot)$ is a symmetric matrix affinely depending on $[x;u]$; sets of this type automatically are convex. In simple words: if you can feed the problem of minimizing a linear function over $X$ to {\tt CVX} \footnote{M.~Grant and S.~Boyd.
  {\em The {\tt CVX} Users’ Guide. {Release} 2.1}, 2014.
  \url{http://web.cvxr.com/cvx/doc/CVX.pdf}}, your $X$ is computationally tractable, and nearly vice versa.
  \subsection{Calculus of compatibility}\label{app:compcalc}
The principal rules of the calculus of compatibility are as follows (verification of the rules is straightforward and is therefore skipped):
\begin{enumerate}
\item {[passing to a subset]} When $\cY'\subset\cY$ are convex compact subsets of $\bR^N$ and a cone $\bY$ is compatible  with $\cY$, the cone is compatible with $\cY'$ as well.
\item {[finite intersection]} Let cones $\bY^j$ be compatible with convex compact sets $\cY_j\subset\bR^N$, $j=1,...,J$. Then the cone
{\small$$
\bY=\cl \{(V,\tau)\in\bS^N_+\times\bR_+:\exists ((V_j,\tau_j)\in\bY^j,j\leq J): V\preceq \sum_jV_j,\sum_j\tau_j\leq \tau\}
$$}\noindent
is compatible with $\cY=\bigcap\limits_j\cY_j$. The closure operation can be skipped when all cones $\bY^j$ are sharp, in which case $\bY$ is sharp as well.
\item {[convex hulls of finite union]} Let cones $\bY^j$ be compatible with convex compact sets $\cY_j\subset\bR^N$, $j=1,...,J$, and let there exist $(V,\tau)$ such that $V\succ0$ and
$$
(V,\tau)\in \bY:=\bigcap\limits_j \bY^j.
$$
Then $\bY$ is compatible with $\cY=\Conv\{\bigcup\limits_j \cY_j\}$ and, in addition, is sharp provided that at least one of  $\bY^j$ is sharp.
\item\label{plhdidp} {[direct product]} Let cones $\bY^j$ be compatible with convex compact sets $\cY_j\subset\bR^{N_j}$, $j=1,...,J$. Then the cone
{\small$$
\bY=\{(V,\tau)\in\bS^{N_1+...+N_J}_+\times \bR_+: \exists (V_j,\tau_j)\in\bY^j: V\preceq\Diag\{V_1,...,V_J\}\ \&\ \tau\geq \sum_j\tau_j\}
$$}\noindent
is compatible with $\cY=\cY_1\times...\times\cY_J$. This cone is sharp, provided that all $\bY^j$ are so.
\item\label{plhdili} {[linear image]} Let cone $\bY$ be compatible with convex compact set  $\cY\subset\bR^N$, let $A$ be a $K\times N$ matrix, and let $\cZ=A\cY$. The cone
$$
\bZ=\cl \{(V,\tau)\in\bS^K_+\times\bR_+: \exists U\succeq A^TVA: (U,\tau)\in\bY\}
$$
is compatible with $\cZ$. The closure operation can be skipped whenever $\bY$ is either sharp, or {\sl complete}, completeness meaning that $(V,\tau)\in\bY$ and $0\preceq V'\preceq V$ imply that $(V',\tau)\in\bY$. The cone $\bZ$ is sharp, provided $\bY$ is so and the rank of $A$ is $K$.
\item {[inverse linear  image]} Let cone $\bY$ be compatible with convex compact set  $\cY\subset\bR^N$, let $A$ be a $N\times K$ matrix with trivial kernel, and let $\cZ=A^{-1}\cY:=\{z\in\bR^K:Az\in\cY\}$. The cone
$$
\bZ=\cl \{(V,\tau)\in\bS^K_+\times\bR_+: \exists U: A^TUA\succeq V \ \&\  (U,\tau)\in\bY\}
$$
is compatible with $\cZ$. The closure operations can be skipped whenever $\bY$ is sharp, in which case $\bZ$ is sharp as well.
\item {[arithmetic summation]} Let cones $\bY^j$ be compatible with convex compact sets $\cY_j\subset\bR^N$, $j=1,...,J$. Then the arithmetic sum
$\cY=\cY_1+...+\cY_J$ of the sets $\cY_j$ can be equipped with compatible cone readily given by the cones $\bY^j$; this cone is sharp, provided all $\bY^j$ are so.\\
Indeed, the arithmetic sum of $\cY_j$ is the linear image of the direct product of $\cY_j$'s under the mapping $[y^1;...;y^J]\mapsto y^1+...+y^J$, and it remains to combine
rules \ref{plhdidp} and \ref{plhdili}; note the cone yielded by rule \ref{plhdidp} is complete, so that when applying rule \ref{plhdili}, the closure operation can be skipped.
\end{enumerate}
\subsection{Proofs}
\subsubsection{Proofs for Section \ref{setscH}}\label{PoissProof}
\paragraph{1$^o$ Sub-Gaussian case.} Note that when $h\in\bR^n$ is deterministic and $\xi\sim \SG(0,\sigma^2I_m)$ we have
\[
\Prob\{|h^T\xi|>1\}\leq 2\exp\left\{-{1\over 2\sigma^2\|h\|_2^2}\right\}.
\]
Indeed, when $h\neq0$ and  $\gamma>0$,
$$
\begin{array}{l}
\Prob\{h^T\xi>1\}\leq \exp\{-\gamma\}\bE\left\{\exp\{\gamma h^T\xi\}\right\}\leq\exp\{{1\over 2}\sigma^2\gamma^2 \|h\|_2^2-\gamma\}.
\end{array}
$$
When minimizing  the resulting bound in $\gamma>0$, we get  $\Prob\{h^T\xi>1\}\leq\exp\left\{-{1\over 2\|h\|_2^2\sigma^2}\right\}$; and the same reasoning
as applied to $-h$ in the role of $h$ results in $\Prob\{h^T\xi<-1\}\leq\exp\left\{-{1\over 2\|h\|_2^2\sigma^2}\right\}$.
Consequently
\[
\pi_G(h):=\underbrace{ \sigma\sqrt{2\ln(2/\delta)}}_{\vartheta_G}\|h\|_2\leq 1\;\;\Rightarrow\;\; \Prob\{|h^T\xi|>1\}\leq\delta,
\]
implying that for $\cH^G_\delta$ as in \rf{eq:gausset} we indeed have
\[
h\in \cH^G_\delta\;\;\Rightarrow\;\;\Prob\{|h^T\xi|>1\}\leq\delta.
\]
\paragraph{2$^o$ Discrete case.} Given $x\in\cX$, setting $\mu=Ax$ and $\eta_k=\zeta_k-\mu$, we get
$$
\omega=Ax+\underbrace{{1\over K}\sum_{k=1}^K\eta_k}_{\xi_x}.
$$
Given $h\in\bR^m$,
$$
h^T\xi_x={1\over K}\sum_k \underbrace{h^T\eta_k}_{\chi_k}.
$$
Random variables $\chi_1,...,\chi_K$ are independent zero mean and clearly satisfy
$$
\bE\left\{\chi_k^2\right\}\leq \sum_i[Ax]_ih_i^2,\,\,|\chi_k|\leq 2\|h\|_\infty.
$$
Applying the Bernstein inequality
we get
\begin{equation}\label{plhd15}
\Prob\{|h^T\xi_x|>1\}=\Prob\left\{|{\sum}_k\chi_k|>K\right\}\leq 2\exp\left\{-{K\over 2\sum_i[Ax]_ih_i^2+{4\over 3}\|h\|_\infty}\right\}.
\end{equation}
Let now
\bse
\pi_D(h)=2\sqrt{\vartheta_D\max_{x\in\cX}\sum_i[Ax]_ih_i^2 +\mbox{$16\over 9$}\vartheta^2_D\|h\|^2_\infty}\;\;\mbox{with}\;\;\vartheta_D={\ln(2/\delta)\over K}.
\ese
After a straightforward computation, we conclude from (\ref{plhd15}) that
\[
\pi_D(h)\leq 1\;\;\Rightarrow\;\; \Prob\{|h^T\xi_x|>1\}\leq \delta,\,\,\forall x\in\cX.
\]
Thus, in the Discrete case we can set
\[
\cH_\delta=\cH_\delta^D:=\{h:\;\pi_D(h)\leq1\}.
\]
\paragraph{Poisson case.} In the Poisson case, for $x\in\cX$, setting $\mu=Ax$, we have
$$
\omega =Ax+\xi_x,\,\xi_x=\omega-\mu.
$$
In this case, for all  $h\in\bR^m$ one has
\begin{equation}\label{EqPoisson}
  \forall t\geq0:
   \Prob\left\{|h^T\xi_x|\geq t\right\}\leq 2\exp\left\{-{t^2\over 2[\sum_ih_i^2\mu_i+\|h\|_\infty t/3]}\right\}
\end{equation}
\begin{quotation}{\small
\noindent Indeed, let $h\in\bR^m$, and let $\omega$ be random vector with independent across $i$ entries $\omega_i\sim\hbox{Poisson}(\mu_i)$.
We have
    \bse
\bE\left\{\exp\{\gamma h^T\omega\}\right\}&=&\prod_i\bE\left\{\gamma h_i\omega_i\right\}=\prod_i\exp\{[\exp\{\gamma h_i\}-1]\mu_i\}\\
    &=&
    \exp\left\{\sum_i[\exp\{\gamma h_i\}-1]\mu_i\right\}.
\ese
 Hence, due to the Markov inequality for $\gamma\geq0$ it holds
    \be
\lefteqn{\Prob\{h^T\omega>h^T\mu+t\}=\Prob\{\gamma h^T\omega>\gamma h^T\mu+\gamma t\}}\nn
&\leq& \bE\left\{\exp\{\gamma h^T\omega\}\right\}\exp\{-\gamma h^T\mu-\gamma t\}
    \leq \exp\left\{\sum_i[\exp\{\gamma h_i\}-\gamma h_i-1]\mu_i-\gamma t\right\}.
    \ee{Eqw10}
Now the standard computation (see, e.g., \cite{massart2007concentration}) shows that
\[
\gamma\leq 3/\|h\|_\infty\;\;\Rightarrow\;\;\left|\e^{\gamma h_i}-\gamma h_i-1\right|\leq {\gamma^2h_i^2\over 2(1-\gamma\|h\|_\infty/3)}
\] which combines with \rf{Eqw10} to imply
\begin{equation}\label{Eqw11}
\ln\left(\Prob\{h^T\omega>h^T\mu+t\}\right)\leq {\gamma^2\sum_ih_i^2\mu_i\over 2(1-\gamma\|h\|_\infty/3)}-\gamma t.
    \end{equation}
For  $\gamma_*={t\over\sum_i h_i^2\mu_i+\|h\|_\infty t/3}\in[0,3/\|h\|_\infty]$ we get:
    $$
    \Prob\left\{h^T\omega>h^T\mu+t\right\}\leq \exp\left\{-{t^2\over 2[\sum_ih_i^2\mu_i+\|h\|_\infty t/3]}\right\}.
    $$
    This inequality combines with the same inequality applied to $-h$ in the role of $h$ to imply (\ref{EqPoisson}). \qed
}\end{quotation}
From \rf{EqPoisson}, we conclude via a straightforward computation that setting
\[
\pi_P(h)=\sqrt{4\vartheta_P \max_{x\in\cX}\sum_i[Ax]_ih_i^2 +\mbox{$16\over 9$} \vartheta^2_P\|h\|^2_\infty}\;\;\mbox{with}\;\;\vartheta_P=\ln(2/\delta),
\]
we ensure that
\[
\pi_P(h)\leq 1\Rightarrow \Prob\{|h^T\xi_x|>1\}\leq \delta,\,\,\forall x\in\cX.
\]
Thus, in the Poisson case we can set
\[
\cH_\delta=\cH_\delta^P:=\{h:\pi_P(h)\leq1\}.
\]
\subsubsection{Proof of Proposition \ref{plhdpropsimp}}\label{proofsimp}
Consider optimization problem specifying $\Psi$ in (\ref{plhdNEq36}). Setting $\theta=r/\rho\geq1$,  let us pass in this problem from variables $v_\ell $ to variables $z_\ell =v_\ell^\rho$, so that
\[\Psi^r=2^r\max_z\left\{\sum_\ell z_\ell^\theta (b_\ell /d_\ell )^r:\;\sum_\ell  z_\ell \leq 1, \;0\leq z_\ell \leq (d_\ell \varsigma_\ell /b_\ell )^\rho\right\}
\leq 2^r\Gamma,
\]
where
\[
\Gamma=\max_z\left\{\sum_\ell z_\ell ^\theta (b_\ell /d_\ell )^r:\;\sum_\ell  z_\ell \leq 1, \;0\leq z_\ell \leq \chi_\ell :=(\vartheta_Gd_\ell /a_\ell )^\rho\right\}\\
\]
(we have used (\ref{plhdtolya11})).
Note that $\Gamma$ is the optimal value in the problem of maximizing a convex (since $\theta\geq1$) function $\sum_\ell z_\ell ^\theta (b_\ell /d_\ell )^r$ over a bounded polyhedral set, so that the maximum is achieved at an extreme point $\bar{z}$ of the feasible set. By the standard characterization of extreme points, the (clearly nonempty) set $I$ of positive entries in $\bar{z}$ is as follows: denoting by $I'$ the set of indexes $\ell
\in I$ such that $\bar{z}_\ell $ is on its upper bound $\bar{z}_\ell =\chi_\ell $, its cardinality $|I'|$ is at least $|I|-1$. Since $\sum_{\ell\in I'}\bar{z}_\ell =
\sum_{\ell\in I'}\chi_\ell \leq1$ and $\chi_\ell $ are nondecreasing in $\ell$ by  (\ref{plhdNEq37}.$b$), we conclude that
$$
\sum_{\ell=1}^{|I'|}\chi_\ell \leq 1,
$$
implying that $|I'|<{\myn}$ when ${\myn}<n$, so that in this case $|I|\leq {\myn}$; and of course $|I|\leq {\myn}$ when ${\myn}=n$. Next, we have
$$
\Gamma=\sum_{\ell\in I}\bar{z}_\ell ^\theta(b_\ell /d_\ell )^r\leq \sum_{\ell\in I}\chi_\ell ^\theta(b_\ell /d_\ell )^r=\sum_{\ell\in I}(\vartheta_Gb_\ell/a_\ell)^r,
$$
and since $b_\ell/a_\ell $ is nonincreasing in $\ell$ and $|I|\leq {\myn}$, the latter quantity is at most $\sum_{\ell=1}^{\myn}(\vartheta_Gb_\ell/a_\ell )^r.$ \qed

\subsubsection{Proof of Lemma \ref{lemThetamu}}\label{lemThetamuproof}
\bigskip\par\noindent\textbf{(i):} When $\pi(\Col_\ell[H]))\leq1$ for all $\ell$ and $\lambda\geq0$, denoting by $[h]^2$ the vector comprised of squares of the entries in $h$, we have
\bse
\phi(\diag(H\Diag\{\lambda\}H^T))&=&\phi(\sum_\ell\lambda_\ell [\Col_\ell[H]]^2)\leq\sum_\ell\lambda_\ell \phi([\Col_\ell[H]]^2)\\
&=&\sum_\ell\lambda_\ell \pi^2(\Col_\ell[H])
\leq\sum_\ell\lambda_\ell,\\
\ese
implying that $(H^T\Diag\{\lambda\}H^T,\;\varkappa\sum_\ell\lambda_\ell)$ belongs to $\bH$.
\bigskip\par\noindent\textbf{(ii):} Let $\Theta,\mu,Q,V$ be as stated in (ii); there is nothing to prove when $\mu=0$, thus assume that $\mu>0$. Let $d=\diag(\Theta)$, so that
\begin{equation}\label{plhdeq60}
d_i=\sum_jQ_{ij}^2,\;\mbox{and}\; \varkappa\phi(d)\leq\mu
\end{equation}
(the second relation is due to $(\Theta,\mu)\in\bH$). (\ref{plhdeq40}) is evident.
We have
$$
[H_\chi]_{ij}=\sqrt{m/\mu}[G_\chi]_{ij},\;\; G_\chi=Q\Diag\{\chi\}V=\left[\sum_{k=1}^m Q_{ik}\chi_kV_{kj}\right]_{i,j}.
$$
We claim that for every $i$ it holds
\begin{equation}\label{change}
\forall \gamma>0:\; \Prob\left\{[G_\chi]_{ij}^2>3\gamma d_i/m\right\}\leq \sqrt{3}\exp\{-\gamma/2\}.
\end{equation}
Indeed, let us fix $i$. There is nothing to prove when $d_i=0$, since in this case $Q_{ij}=0$ for all $j$ and therefore $[G_\chi]_{ij}\equiv 0$. When $d_i>0$, by homogeneity in $Q$, it suffices to verify (\ref{change}) when $d_i/m=1/3$. Assuming that this is the case,
let $\eta\sim\cN(0,1)$ be independent of $\chi$. We have
\bse
\bE_{\eta}\left\{\bE_{\chi}\{\exp\{\eta [G_\chi]_{ij}\}\}\right\}&=&\bE_{\eta}\left\{\prod_{k}\cosh(\eta Q_{ik}V_{kj})\right\}
\leq \bE_\eta\left\{\prod_k\exp\{\half\eta^2Q_{ik}^2V_{kj}^2\}\right\}\\
&=&\bE_\eta\big\{\exp\{\half\eta^2\underbrace{{\sum}_kQ_{ik}^2V_{kj}^2}_{\leq 2d_i/m}\}\big\}\leq\bE_\eta\left\{\eta^2d_i/m\right\}=\bE_\eta\left\{\exp\{\eta^2/3\}\right\}=\sqrt{3}.
\ese
On the other hand,
$$
\bE_{\chi}\left\{\bE_{\eta}\{\exp\{\eta [G_\chi]_{ij}\}\}\right\}=\bE_\chi\left\{\exp\{\half [G_\chi]_{ij}^2\}\right\},
$$
implying that
$$
\bE_\chi\left\{\exp\{\half[G_\chi]_{ij}^2\}\right\}\leq\sqrt{3}.
$$
Therefore in the case of $d_i/m=1/3$ for all $s>0$ it holds
$$
\Prob\{\chi: [G_\chi]_{ij}^2>s\}\leq \sqrt{3}\exp\{-s/2\},
$$
and (\ref{change}) follows. Recalling the relation between $H$ and $G$, we get from (\ref{change}) that for all $\gamma>0$
$$
\Prob\{\chi: [H_\chi]_{ij}^2>3\gamma d_i/\mu\}\leq \sqrt{3}\exp\{-\gamma/2\}.
$$
Therefore, with $\varkappa$ given by (\ref{plhdeq50}),
the probability of the event
$$
\forall i,j: [H_\chi]_{ij}^2\leq \varkappa{d_i\over \mu}
$$
is at least 1/2. Let this event take place; in this case we have $[\Col_\ell[H_\chi]]^2\leq \varkappa d/\mu$, whence, by definition of the norm $\pi(\cdot)$,
$\pi^2(\Col_\ell[H_\chi])\leq \varkappa \phi(d)/\mu\leq1$ (see the inequality in (\ref{plhdeq60})). Thus, the probability of the event
 (\ref{plhdeq41}) is at least 1/2. \qed
\subsubsection{
Verification of (\ref{plhdEq20}), (\ref{plhdEq20a})}\label{AppSpecEll}
Let us verify that (\ref{plhdEq20}) is a closed convex cone compatible with the spectratope
$\cX=M\cY$, with $\cY$ given by (\ref{ellspectr}.$b$); verification of the ``ellitopic'' version of
this claim can be obtained from what follows by straightforward  simplifications. The only non-evident fact in the claim  is
that whenever $(V,\tau)\in\bX$, we have
\begin{equation}\label{wehave3432}
x^TVx\leq\tau\,\,\forall x\in\cX.
\end{equation}
To justify (\ref{wehave3432}), let $(V,\tau)\in\bX$, so that
\begin{equation}\label{wehave3433}
M^TVM\preceq \sum_\ell \cR_\ell^*[\Lambda_\ell]\;\mbox{and}\;\phi_\cR(\lambda[\Lambda])\leq\tau
\end{equation}
for properly selected $\Lambda_\ell\in\bS^{d_\ell}$, and let $x\in\cX$, so that $x=My$ for some $y$ satisfying the relations
\begin{equation}\label{wehave3434}
\sum_{i,j}y_iy_j R^{\ell i}R^{\ell j}=\aic{\R}{R}_\ell^2[y]\preceq r_\ell I_{d_\ell},\;\ell\leq L
\end{equation}
for some properly selected $r\in\cR$.
Taking into account that $R^{\ell i}$ are symmetric and $\Lambda_\ell\succeq0$, (\ref{wehave3434}) implies that
$$
y^T\cR^*_\ell[\Lambda_\ell]y=\sum_{i,j}y_iy_j\Tr(R^{\ell i}\Lambda_\ell R^{\ell j})=\Tr(R_\ell[y]\Lambda_\ell R_\ell[y])=
\Tr(\Lambda_\ell R_\ell^2[y])\leq r_\ell\Tr(\lambda_\ell),
$$
which combines with  the first relation in (\ref{wehave3433}) to imply that
$$
x^TVx=y^T[M^TVM]y\leq \sum_\ell y^T\cR^*_\ell[\Lambda_\ell]y\leq \sum_\ell r_\ell\Tr(\lambda_\ell)\leq \phi_\cR(\lambda[\Lambda]),
$$
where the concluding inequality is due to $r\in\cR$. Invoking the second relation in (\ref{wehave3433}) and recalling that $x\in\cX$ is arbitrary,
we arrive at (\ref{wehave3432}). \qed

  \subsubsection{Verification of (\ref{plhdEq204b})}\label{plhdeasy}
Given $s\in[2,\infty]$ and setting $\bar{s}=s/2$, $s_*={s\over s-1}$, $\bar{s}_*={\bar{s}\over\bar{s}-1}$, we want to prove that
\[\begin{array}{l}
\left\{(V,\tau)\in\bS^N_+\times\bR_+:\;\exists (W\in\bS^N,w\in\bR^N_+): \;V\preceq W+\Diag\{w\}\; \&\; \|W\|_{s_*}+\|w\|_{\bar{s}_*}\leq\tau\right\}\\
\quad=~\left\{(V,\tau)\in\bS^N_+\times\bR_+:\;\exists w\in\bR^N_+:\;V\preceq\Diag\{w\},
\|w\|_{\bar{s}_*}\leq\tau\right\}.
\end{array}
\]
To this end it suffices to check that whenever $W\in\bS^N$ there exists $w\in\bR^N$ satisfying
$$
W\preceq\Diag\{w\},\;\;\|w\|_{\bar{s}_*}\leq \|W\|_{s_*}.
$$
The latter claim is nothing but the claim that whenever $W\in\bS^N$, and $\|W\|_{s_*}\leq1$, the conic optimization problem
\begin{equation}\label{Eq1000}
\Opt=\min_{t,w}\{t:t\geq \|w\|_{\bar{s}_*},\Diag\{w\}\succeq W\}
\end{equation}
is solvable (which is evident) with optimal value $\leq 1$. To see that the latter indeed is the case, note that the problem clearly is strictly feasible, whence its optimal value is the same as the optimal value in the conic problem
$$
\begin{array}{c}
\Opt=\max_P\left\{\Tr(PW):P\succeq0,\, \|\diag(P)\|_{\bar{s}_*/(\bar{s}_*-1)}\leq 1\right\}\\
\left[\diag(P)=[P_{11};P_{22};...;P_{NN}]\right]\\
\end{array}
$$
dual to (\ref{Eq1000}). Since
\[
\Tr(PW)\leq \|P\|_{s_*/(s_*-1)}\|W\|_{s_*}\leq \|P\|_{s_*/(s_*-1)},
\]
when recalling what $s_*$ and $\bar{s}_*$ are, our task boils down to verifying that whenever a matrix $P\succeq0$ satisfies $\|\diag(P)\|_{s/2}\leq1$, one has also $\|P\|_s\leq1$. This is immediate: since $P$ is positive semidefinite, we have $|P_{ij}|\leq P_{ii}^{1/2}P_{jj}^{1/2}$, whence, assuming $s<\infty$,
$$
\|P\|_s^s=\sum_{i,j}|P_{ij}|^s\leq\sum_{i,j}P_{ii}^{s/2}P_{jj}^{s/2}=\left(\sum_iP_{ii}^{s/2}\right)^2\leq 1.
$$
When $s=\infty$, the same argument leads to
$$\|P\|_\infty=
\max_{i,j}|P_{ij}|=\max_{i} |P_{ii}|=\|\diag(P)\|_\infty. \eqno{\mbox{\qed}}
$$
\subsubsection{Proof of Proposition \ref{pr:near}}\label{proofnear}
{\bf 1$^o$.}
Let us consider the optimization problem, as defined in \cite[relation (26)]{JudNem2018} (where one should set $\cQ=\sigma^2 I_m$) which under the circumstances is responsible for building a nearly optimal {\sl linear} estimate of $w=Bx$, namely,
\begin{equation}\label{plhdeq1000}
\begin{array}{ll}
\Opt_*&=\min\limits_{\Theta,H,\Lambda,\Upsilon',\Upsilon''}\bigg\{\phi_\cT(\lambda[\Lambda])
+\phi_\cR(\lambda[\Upsilon'])+\phi_\cR(\lambda[\Upsilon''])+\sigma^2\Tr(\Theta):\\
&\begin{array}{l}\Lambda=\{\Lambda_k\succeq0,k\leq K\},\Upsilon'=\{\Upsilon_\ell^\prime\succeq 0,\ell\leq L\},\Upsilon''=\{\Upsilon_\ell^{\prime\prime}\succeq0,\ell\leq L\},\\
\left[\begin{array}{c|c}\sum_\ell\cS_\ell^*[\Upsilon_\ell^\prime]&\half M^T[B-H^TA]\cr\hline
\half [B-H^TA]^TM&\sum_k\cR_k^*[\Lambda_k]\cr\end{array}\right]\succeq0,\;
\left[\begin{array}{c|c}\sum_\ell\cS_\ell^*[\Upsilon_\ell^{\prime\prime}]&\half M^TH^T\cr\hline\half HM&\Theta\cr\end{array}\right]\succeq0\\
\end{array}\Bigg\}
\end{array}
\end{equation}
Let us show that the optimal value $\Opt$ of \rf{plhdreadsini} satisfies
\begin{equation}\label{plhdconclude}
\Opt\leq 2\kappa\Opt_*=2\sqrt{2\ln(2m/\epsilon)}\Opt_*.
\end{equation}
To this end, observe that the matrices
\[
Q:=\left[\begin{array}{c|c}U&{1\over 2}B\cr\hline {1\over 2}B^T&A^T\Theta A+\sum_k \cR_k^*[\Lambda_k]\cr\end{array}\right]
\] and
$$
\left[\begin{array}{c|c}M^TUM&{1\over 2}M^TB\cr\hline {1\over 2}B^TM&A^T\Theta A+\sum_k \cR_k^*[\Lambda_k]\cr\end{array}\right]=\left[\begin{array}{c|c}M^T&\cr\hline&I_n\cr\end{array}\right]Q\left[\begin{array}{c|c}M&\cr\hline&I_n\cr\end{array}\right]
$$
simultaneously are/are not positive semidefinite due to the fact that the image space of $M$ contains the full-dimensional set $\cB_*$ and thus is the entire $\bR^\nu$, so that the image space of
$\left[\begin{array}{c|c}M&\cr\hline&I_n\cr\end{array}\right]$ is the entire $\bR^\nu\times\bR^n$. Therefore
\begin{equation}\label{plhdreads}
\begin{array}{rcl}
\Opt&=&\min\limits_{\Theta,U,\Lambda,\Upsilon}\bigg\{2\left[\phi_\cR(\lambda[\Upsilon])+\phi_{\cT}(\lambda[\Lambda])+\sigma^2\kappa^2\Tr(\Theta)\right]:\\
&&\begin{array}{l}\Theta\succeq0,U\succeq0,\Lambda=\{\Lambda_k\succeq0,k\leq K\},\Upsilon=\{\Upsilon_\ell\succeq0,\ell\leq L\},\\
\left[\begin{array}{c|c}M^TUM&{1\over 2}M^TB\cr\hline {1\over 2}B^TM&A^T\Theta A+\sum_k \cR_k^*[\Lambda_k]\cr\end{array}\right]\succeq0,\,\,
M^TUM\preceq \sum_\ell\cS^*_\ell[\Upsilon_\ell]\\
\end{array}\Bigg\}
\end{array}
\end{equation}
Further, note that if a collection $\Theta,U,\{\Lambda_k\},\{\Upsilon_\ell\}$ is a feasible solution to the latter problem and $\theta>0$, the scaled collection $\theta\Theta,\theta^{-1}U,\{\theta\Lambda_k\},
\{\theta^{-1}\Upsilon_\ell\}$ is also a feasible solution. When optimizing with respect to the scaling, we get
\be
\Opt&=&\inf\limits_{\Theta,U,\Lambda,\Upsilon}
\bigg\{4\sqrt{\phi_\cR(\lambda[\Upsilon])\left[\phi_\cT(\lambda[\Lambda]
+\sigma^2\kappa^2\Tr(\Theta)\right]}:\nn
&&\begin{array}{l}\Theta\succeq0,U\succeq0,\Lambda=\{\Lambda_k\succeq0,k\leq K\},\Upsilon=\{\Upsilon_\ell\succeq0,\ell\leq L\},\\
\left[\begin{array}{c|c}M^TUM&{1\over 2}M^TB\cr\hline {1\over 2}B^TM&A^T\Theta A+\sum_k \cR_k^*[\Lambda_k]\cr\end{array}\right]\succeq0,\,\,
M^TUM\preceq \sum_\ell\cS^*_\ell[\Upsilon_\ell]\\
\end{array}\bigg\}
\nn&\leq&2\kappa\Opt_+,
\ee{eq:64p}
where
\begin{eqnarray}\label{eq:69}
\Opt_+&=&\inf\limits_{\Theta,U,\Lambda,\Upsilon}
\bigg\{2\sqrt{\phi_\cR(\lambda[\Upsilon])\left[\phi_{\cT}(\lambda[\Lambda])+\sigma^2\Tr(\Theta)\right]}:\nn
&&\begin{array}{l}\Theta\succeq0,U\succeq0,\Lambda=\{\Lambda_k\succeq0,k\leq K\},\Upsilon=\{\Upsilon_\ell\succeq0,\ell\leq L\},\\
\left[\begin{array}{c|c}M^TUM&{1\over 2}M^TB\cr\hline {1\over 2}B^TM&A^T\Theta A+\sum_k \cR_k^*[\Lambda_k]\cr\end{array}\right]\succeq0,\,\,
M^TUM\preceq \sum_\ell\cS^*_\ell[\Upsilon_\ell]\\
\end{array}\bigg\}\\
&&\hbox{[note that $\kappa>1$]}\nonumber
\end{eqnarray}
On the other hand, when strengthening the constraint $\Lambda_k\succeq0$ of \rf{plhdeq1000} to $\Lambda_k\succ 0$, we still have
\begin{equation}\label{plhdeq100}
\begin{array}{rcl}
\Opt_*&=\inf\limits_{\Theta,H,\Lambda,\Upsilon',\Upsilon''}\bigg\{\phi_\cT(\lambda[\Lambda])+\phi_\cR(\lambda[\Upsilon'])+\phi_\cR(\lambda[\Upsilon''])+\sigma^2\Tr(\Theta):\\
&\begin{array}{l}\Lambda=\{\Lambda_k\succ0,k\leq K\},\Upsilon'=\{\Upsilon_\ell^\prime\succeq0,\ell\leq L\},\Upsilon''=\{\Upsilon_\ell^{\prime\prime}\succeq0,\ell\leq L\},\\
\left[\begin{array}{c|c}\sum_\ell\cS_\ell^*[\Upsilon_\ell^\prime]&\half M^T[B-H^TA]\cr\hline
\half [B-H^TA]^TM&\sum_k\cR_k^*[\Lambda_k]\cr\end{array}\right]\succeq0,\;
\left[\begin{array}{c|c}\sum_\ell\cS_\ell^*[\Upsilon_\ell^{\prime\prime}]&\half M^TH^T\cr\hline\half HM&\Theta\cr\end{array}\right]\succeq0\\
\end{array}\Bigg\}.
\end{array}
\end{equation}
Now let $\Theta,H,\Lambda,\Upsilon',\Upsilon''$ be a feasible solution to the latter problem. By the second semidefinite constraint in
(\ref{plhdeq100}) we have
\bse
\left[\begin{array}{c|c}\sum_\ell\cS_\ell^*[\Upsilon_\ell^{\prime\prime}]&\half M^TH^TA\cr\hline\half
A^THM&A^T\Theta A\cr\end{array}\right]&=&\left[\begin{array}{c|c}I&\cr\hline&A\cr\end{array}\right]^T\left[\begin{array}{c|c}
\sum_\ell\cS_\ell^*[\Upsilon_\ell^{\prime\prime}]&\half M^TH^T\cr\hline
\half HM&\Theta\cr\end{array}\right]\left[\begin{array}{c|c}I&\cr\hline&A\cr\end{array}\right]
\succeq0
\ese
which combines with the first semidefinite constraint in (\ref{plhdeq100}) to imply that
$$
\left[\begin{array}{c|c}\sum_\ell\cS_\ell^*[\Upsilon_\ell^\prime+\Upsilon_\ell^{\prime\prime}]&\half
M^TB\cr\hline\half B^TM&A^T\Theta A+\sum_k\cR_k^*[\Lambda_k]\cr\end{array}\right]\succeq0.
$$
Next, by the Schur Complement Lemma (which is applicable due to
\[A^T\Theta A+\sum_k\cR_k^*[\Lambda_k]\succeq\sum_k\cR_k^*[\Lambda_k]\succ0,
\] where the concluding $\succ$ is due to \cite[Lemma 5.1]{JudNem2018} combined with $\Lambda_k\succ0$), this relation implies that for
$$\Upsilon_\ell=\Upsilon_\ell^\prime+\Upsilon_\ell^{\prime\prime},
$$
we have
$$
\sum_\ell\cS_\ell^*[\Upsilon_\ell]\succeq M^T\underbrace{\left[\four B[A^T\Theta A+\sum_k\cR_k^*[\Lambda_k]]^{-1}B^T\right]}_{U}M.
$$
Using the Schur Complement Lemma again, for the just defined $U\succeq0$ we obtain
$$
\left[\begin{array}{c|c}M^TUM&{1\over 2}M^TB\cr\hline{1\over 2}B^TM&A^T\Theta A+\sum_k\cR_k^*[\Lambda_k]\cr\end{array}\right]\succeq0,
$$
and in addition, by the definition of $U$,
$$
M^TUM\preceq \sum_\ell\cS_\ell^*[\Upsilon_\ell].
$$
 We conclude that $$
(\Theta,U,\Lambda,\Upsilon:=\{\Upsilon_\ell=\Upsilon_\ell^\prime+\Upsilon_\ell^{\prime\prime},\ell\leq L\})$$
is a feasible solution to optimization problem (\ref{eq:69}) specifying $\Opt_+$. The value of the objective of the latter problem at this feasible solution is
\bse
2\sqrt{\phi_\cR(\lambda[\Upsilon']+\lambda[\Upsilon''])\left[\phi_\cT(\lambda[\Lambda]) +\sigma^2\Tr(\Theta)\right]}&
\leq&
\phi_\cR(\lambda[\Upsilon']+\lambda[\Upsilon''])+\phi_\cT(\lambda[\Lambda]) +\sigma^2\Tr(\Theta)\\
&\leq& \phi_\cR(\lambda[\Upsilon'])+\phi_\cR(\lambda[\Upsilon''])+\phi_\cT(\lambda[\Lambda])
+\sigma^2\Tr(\Theta),
\ese
the concluding quantity in the chain being the value of the objective of problem (\ref{plhdeq100}) at the feasible solution
$\Theta,H,\Lambda,\Upsilon',\Upsilon''$ to this problem. Since the resulting inequality holds true for every feasible solution to (\ref{plhdeq100}), we
conclude that $\Opt_+\leq\Opt_*$, and we arrive at \rf{plhdconclude} due to \rf{eq:64p}.
\paragraph{2$^o$.} Now, from Proposition 3.3 in the latest arXiv version of \cite[Section 5.7]{JudNem2018}, we conclude that $\Opt_*$ is within a logarithmic factor  of the minimax optimal $({1\over 8},\|\cdot\|)$-risk corresponding to the case of Gaussian noise $\xi_x\sim\cN(0,\sigma^2I_m)$ for all $x$:
\[
\Opt_*\leq \theta_*\RiskOpt_{1/8},
\]
where \[\theta_*=8\sqrt{\left(2 \ln F+10\ln 2\right)
\left(2 \ln D+10\ln 2\right)},\quad F=\sum_\ell f_\ell,\;D=\sum_k d_k.
\]
Since the minimax optimal $(\epsilon,\|\cdot\|)$-risk clearly only grows when $\epsilon$ decreases, we conclude that for $\epsilon\leq1/8$ a feasible near optimal solution to \rf{plhdreadsini} is minimax optimal within the factor $2\theta^*\kappa$.\qed

\end{document}